\documentclass{article}
\usepackage{amsmath,amssymb}
\usepackage{cite}

\begin{document}

\begin{center}
{\bf Cramer's rules  for the solution to the two-sided restricted quaternion matrix equation.}\end{center}
\begin{center}{ \textbf{\emph{Ivan I. Kyrchei} }\footnote{kyrchei@online.ua,\\ Pidstrygach Institute for Applied Problems of Mechanics and Mathematics of NAS of Ukraine,
 Ukraine} }\end{center}

\begin{abstract}
Weighted  singular value decomposition (WSVD) of a quaternion matrix and with its help determinantal representations  of the quaternion weighted Moore-Penrose inverse  have been derived recently by the author.
 In this paper, using these determinantal representations, explicit determinantal representation formulas  for the solution of the  restricted quaternion matrix equations, ${\bf A}{\bf X}{\bf B}={\bf D}$, and consequently, ${\bf A}{\bf X}={\bf D}$ and ${\bf X}{\bf B}={\bf D}$ are obtained within the framework of the theory of  column-row determinants.  We   consider all possible cases  depending on weighted matrices.  \end{abstract}



\section{Introduction}
\newtheorem{corollary}{Corollary}[section]
\newtheorem{theorem}{Theorem}[section]
\newtheorem{lemma}{Lemma}[section]
\newtheorem{definition}{Definition}[section]
\newtheorem{remark}{Remark}[section]
\newcommand{\rank}{\mathop{\rm rank}\nolimits}
\newtheorem{proposition}{Proposition}[section]

Let ${\rm
{\mathbb{R}}}$ and ${\rm
{\mathbb{C}}}$ be the real and complex number fields, respectively.
Throughout the paper, we denote  the set of all $m\times n$ matrices over the quaternion
skew field
\[{\rm {\mathbb{H}}}=\{a_{0}+a_{1}i+a_{2}j+a_{3}k\,
|\,i^{2}=j^{2}=k^{2}=-1,\, a_{0}, a_{1}, a_{2}, a_{3}\in{\rm
{\mathbb{R}}}\}\]
by ${\rm {\mathbb{H}}}^{m\times n}$, and by ${\rm {\mathbb{H}}}^{m\times n}_{r}$ the set of all $m\times n$ matrices over $\mathbb{H}$ with a rank $r$. Let ${\rm M}\left( {n,{\rm {\mathbb{H}}}} \right)$ be the
ring of $n\times n$ quaternion matrices and $ {\bf I}$ be the identity matrix with the appropriate size. For ${\rm {\bf A}}
 \in {\rm {\mathbb{H}}}^{n\times m}$, we denote by ${\rm {\bf A}}^{ *}$, $\rank {\bf A}$ the conjugate transpose (Hermitian adjoint) matrix and the rank
of ${\rm {\bf A}}$.
 The matrix ${\rm {\bf A}} = \left( {a_{ij}}  \right) \in {\rm
{\mathbb{H}}}^{n\times n}$ is Hermitian if ${\rm {\bf
A}}^{ *}  = {\rm {\bf A}}$.

The  definitions of the Moore-Penrose inverse \cite{pen} and the weighted Moore-Penrose inverse \cite{ben} can be  extended  to quaternion matrices as follows.

The Moore-Penrose inverse of ${\bf A}\in{\rm {\mathbb{H}}}^{m\times n}$, denoted by ${\bf A}^{\dagger}$, is the unique matrix ${\bf X}\in{\rm {\mathbb{H}}}^{n\times m}$ satisfying the following equations \cite{pen},
 \begin{gather}\label{eq1:MP}  {\rm {\bf A}}{\bf X}
{\rm {\bf A}} = {\rm {\bf A}}; \\
                                 \label{eq2:MP}  {\bf X} {\rm {\bf
A}}{\bf X}  = {\bf X};\\
                                  \left( {\rm {\bf A}}{\bf X} \right)^{ *}  = {\rm
{\bf A}}{\bf X}; \\
                                \left( {{\bf X} {\rm {\bf A}}} \right)^{ *}  ={\bf X} {\rm {\bf A}}. \end{gather}

Let Hermitian positive definite matrices ${\bf M}$ and ${\bf N}$ of order $m$ and $n$, respectively, be given. For
 ${\bf A}\in {\mathbb H}^{m\times n}$, \textbf{the weighted Moore-Penrose inverse} of ${\bf A}$ is the unique
solution ${\bf X}={\bf A}^{\dag}_{M,N}$ of the matrix equations (\ref{eq1:MP}) and (\ref{eq2:MP}) and the following
equations in ${\bf X}$ \cite{pr}:
\[(3M)\,\,({\bf M}{\rm {\bf A}}{\bf X})^{*}  = {\bf M}{\rm
{\bf A}}{\bf X};\,\,(4N)\,\,({\bf N}{\bf X}{\rm {\bf A}})^{*}  = {\bf N}{\bf X}{\rm
{\bf A}}.\]
In particular, when ${\bf M}={\bf I}_m$ and ${\bf N}={\bf I}_n$, the matrix ${\bf X}$ satisfying the  equations (\ref{eq1:MP}), (\ref{eq2:MP}), (3M), (4N) is the  Moore-Penrose inverse ${\bf A}^{\dag}$.

A basic method for finding the Moore-Penrose inverse is based on the  singular value decomposition (SVD). It is also available for quaternion matrices, (see, e.g. \cite{zha,ky_math_sci}). The  weighted Moore-Penrose inverse ${\bf A}^{\dag}_{M,N}\in {\rm {\mathbb{C}}}^{m\times n} $ (over complex or real fields) has the explicit expressing by the weighted singular value decomposition (WSVD) that at first has been obtained in \cite{loan} by Cholesky factorization. In \cite{galba}, WSVD of
real matrices with singular weights has been derived using
weighted orthogonal matrices and weighted pseudoorthogonal matrices.

Recently, by the author, WSVD has been  expanded to quaternion matrices.
\begin{theorem}\cite{kyr_wmpi}\label{tm_weig_A+} Let ${\bf A}\in {\mathbb{H}}^{m\times n}_{r}$,  ${\bf M}$ and ${\bf N}$ be positive definite matrices  of order $m$ and $n$, respectively. Denote ${\bf A}^{\sharp}={\bf N}^{-1}{\bf A}^{*}{\bf M}$. There exist ${\bf U}\in {\mathbb{H}}^{m\times m}$, ${\bf V}\in {\mathbb{H}}^{n\times n}$ satisfying ${\bf U}^{*}{\bf M}{\bf U}={\bf I}_{m}$ and ${\bf V}^{*}{\bf N}^{-1}{\bf V}={\bf I}_{n}$ such that
${\bf A}={\bf U}{\bf D}
{\bf V}^{*}$, where ${\bf D}=\left(
                 \begin{array}{cc}
                   {\bf \Sigma} & {\bf 0} \\
                   {\bf 0} & {\bf 0} \\
                 \end{array}
               \right)$.
Then the weighted Moore-Penrose inverse ${\bf A}^{\dag}_{M,N}$ can be represented
\begin{equation}\label{eq:wsvd_A+}
{\bf A}^{\dag}_{M,N}={\bf N}^{-1}{\bf V}\left(
                 \begin{array}{cc}
                   {\bf \Sigma}^{-1} & {\bf 0} \\
                   {\bf 0} & {\bf 0} \\
                 \end{array}
               \right)
{\bf U}^{*}{\bf M},
\end{equation}
where ${\bf \Sigma}=diag(\sigma_{1}, \sigma_{2},...,\sigma_{r})$, $\sigma_{1}\geq \sigma_{2}\geq...\geq \sigma_{r}>0$ and $\sigma^{2}_{i}$ is the nonzero eigenvalues of ${\bf A}^{\sharp}{\bf A}$ or  ${\bf A}{\bf A}^{\sharp}$, which coincide.
\end{theorem}

By using WSVD, within the framework of the theory of   column-row determinants,  limit and determinantal representations  of the quaternion weighted Moore-Penrose inverse   has been derived ibidem as well.

But why determinantal representations of generalized inverses are so important? When we return to the usual inverse, its determinantal representation  is the matrix with cofactors in entries that gives direct method of its  finding  and makes it  applicable in Cramer's rule for  systems of linear equations. The same be wanted for generalized inverses.   But there is not so unambiguous even  for complex or real matrices. Therefore, there are various determinantal representations of generalized inverses because searches  of their explicit more applicable expressions are continuing (see, e.g. \cite{st1,liu2,liu1,ky_lma1,ky_nova}).

The understanding of the problem for determinantal representing of generalized inverses as well as solutions   and generalized inverse solutions of  quaternion matrix equations,  only now  begins to be decided due to the theory of column-row determinants introduced in  \cite{kyr2,ky_th}.

Song at al. \cite{song1,song6} have studied the weighted  Moore-Penrose inverse over the quaternion skew field and obtained its determinantal representation  within the framework of the theory of  column-row determinants as well. But  WSVD of
quaternion matrices has not been considered and for obtaining  a determinantal representation there was used auxiliary matrices which different from ${\bf A}$, and weights ${\bf M}$ and ${\bf N}$. Despite this in \cite{song6}, Cramer's rule of the quaternion restricted matrix equation ${\bf A}{\bf X}{\bf B}={\bf D}$ has been derived with the help obtained determinantal representations of the weighted  Moore-Penrose inverse.

The main goals of the paper are   obtaining   Cramer's rule for   the quaternion restricted matrix equation ${\bf A}{\bf X}{\bf B}={\bf D}$, and consequently, ${\bf A}{\bf X}={\bf D}$ and ${\bf X}{\bf B}={\bf D}$ using the determinantal representations of the weighted  Moore-Penrose inverse obtained by WSVD in \cite{kyr_wmpi}. We consider all possible cases with respect to weights of ${\bf A}$ and ${\bf B}$.

It need to note that currently the theory of column-row determinants of quaternion matrices  is active  developing.
Within the framework of  column-row determinants, determinantal representations of various kind of generalized inverses,  (generalized inverses) solutions    of quaternion matrix equations recently have been derived as by the author (see, e.g.\cite{kyr5,kyr6,kyr7,kyr8,kyr9}) so by other researchers (see, e.g.\cite{ky_nova1,song3,song4,song5}).

In this chapter we shall adopt the following notation.

Let $\alpha : = \left\{
{\alpha _{1} ,\ldots ,\alpha _{k}} \right\} \subseteq {\left\{
{1,\ldots ,m} \right\}}$ and $\beta : = \left\{ {\beta _{1}
,\ldots ,\beta _{k}} \right\} \subseteq {\left\{ {1,\ldots ,n}
\right\}}$ be subsets of the order $1 \le k \le \min {\left\{
{m,n} \right\}}$. By ${\rm {\bf A}}_{\beta} ^{\alpha} $ denote the
submatrix of ${\rm {\bf A}}$ determined by the rows indexed by
$\alpha$, and the columns indexed by $\beta$. Then, ${\rm {\bf
A}}{\kern 1pt}_{\alpha} ^{\alpha}$ denotes a principal submatrix
determined by the rows and columns indexed by $\alpha$.
 If ${\rm {\bf A}} \in {\rm
M}\left( {n,{\rm {\mathbb{H}}}} \right)$ is Hermitian, then by
${\left| {{\rm {\bf A}}_{\alpha} ^{\alpha} } \right|}$ denote the
corresponding principal minor of $\det {\rm {\bf A}}$, since ${\rm {\bf
A}}{\kern 1pt}_{\alpha} ^{\alpha}$ is Hermitian as well.
 For $1 \leq k\leq n$, denote by $\textsl{L}_{ k,
n}: = {\left\{ {\,\alpha :\alpha = \left( {\alpha _{1} ,\ldots
,\alpha _{k}} \right),\,{\kern 1pt} 1 \le \alpha _{1} \le \ldots
\le \alpha _{k} \le n} \right\}}$ the collection of strictly
increasing sequences of $k$ integers chosen from $\left\{
{1,\ldots ,n} \right\}$. For fixed $i \in \alpha $ and $j \in
\beta $, let
\[I_{r,\,m} {\left\{ {i} \right\}}: = {\left\{ {\,\alpha :\alpha \in
L_{r,m} ,i \in \alpha}  \right\}}{\rm ,} \quad J_{r,\,n} {\left\{
{j} \right\}}: = {\left\{ {\,\beta :\beta \in L_{r,n} ,j \in
\beta}  \right\}}.\]

 The paper is organized as follows. We start with some basic concepts and results from the theories of   row-column determinants and of quaternion matrices in Section 2.
 Cramer's rules for the quaternionic restricted matrix equation ${\bf A}{\bf X}{\bf B}={\bf D}$, and consequently, ${\bf A}{\bf X}={\bf D}$ and ${\bf X}{\bf B}={\bf D}$ are derived in Section 3. All possible cases are considered in the three subsections of Section 3.
In Section 4, we give
 numerical an example to illustrate the main results.

\section{Preliminaries}
\subsection{Elements of the theory of  column-row determinants and quaternion inverse matrices}
 For a quadratic matrix ${\rm {\bf A}}=(a_{ij}) \in {\rm
M}\left( {n,{\mathbb{H}}} \right)$ can be define $n$ row determinants and $n$ column determinants as follows.

Suppose $S_{n}$ is the symmetric group on the set $I_{n}=\{1,\ldots,n\}$.
\begin{definition}\cite{kyr2}
 The $i$th row determinant of ${\rm {\bf A}}=(a_{ij}) \in {\rm
M}\left( {n,{\mathbb{H}}} \right)$ is defined  for all $i = 1,\ldots,n $
by putting
 \begin{gather*}{\rm{rdet}}_{ i} {\rm {\bf A}} =
{\sum\limits_{\sigma \in S_{n}} {\left( { - 1} \right)^{n - r}({a_{i{\kern
1pt} i_{k_{1}}} } {a_{i_{k_{1}}   i_{k_{1} + 1}}} \ldots } } {a_{i_{k_{1}
+ l_{1}}
 i}} ) \ldots  ({a_{i_{k_{r}}  i_{k_{r} + 1}}}
\ldots  {a_{i_{k_{r} + l_{r}}  i_{k_{r}} }}),\\
\sigma = \left(
{i\,i_{k_{1}}  i_{k_{1} + 1} \ldots i_{k_{1} + l_{1}} } \right)\left(
{i_{k_{2}}  i_{k_{2} + 1} \ldots i_{k_{2} + l_{2}} } \right)\ldots \left(
{i_{k_{r}}  i_{k_{r} + 1} \ldots i_{k_{r} + l_{r}} } \right),\end{gather*}
with
conditions $i_{k_{2}} < i_{k_{3}}  < \ldots < i_{k_{r}}$ and $i_{k_{t}}  <
i_{k_{t} + s} $ for $t = {2,\ldots,r} $ and $s ={1,\ldots,l_{t}} $.
\end{definition}
\begin{definition}\cite{kyr2}
The $j$th column determinant
 of ${\rm {\bf
A}}=(a_{ij}) \in {\rm M}\left( {n,{\mathbb{H}}} \right)$ is defined for
all $j =1,\ldots,n $ by putting
 \begin{gather*}{\rm{cdet}} _{{j}}\, {\rm {\bf A}} =
{{\sum\limits_{\tau \in S_{n}} {\left( { - 1} \right)^{n - r}(a_{j_{k_{r}}
j_{k_{r} + l_{r}} } \ldots a_{j_{k_{r} + 1} i_{k_{r}} } ) \ldots } }(a_{j\,
j_{k_{1} + l_{1}} }  \ldots  a_{ j_{k_{1} + 1} j_{k_{1}} }a_{j_{k_{1}}
j})},\\
\tau =
\left( {j_{k_{r} + l_{r}}  \ldots j_{k_{r} + 1} j_{k_{r}} } \right)\ldots
\left( {j_{k_{2} + l_{2}}  \ldots j_{k_{2} + 1} j_{k_{2}} } \right){\kern
1pt} \left( {j_{k_{1} + l_{1}}  \ldots j_{k_{1} + 1} j_{k_{1} } j}
\right), \end{gather*}
\noindent with conditions, $j_{k_{2}}  < j_{k_{3}}  < \ldots <
j_{k_{r}} $ and $j_{k_{t}}  < j_{k_{t} + s} $ for  $t = {2,\ldots,r} $
and $s = {1,\ldots,l_{t}}  $.
\end{definition}

Suppose ${\rm {\bf A}}_{}^{i{\kern 1pt} j} $ denotes the submatrix of
${\rm {\bf A}}$ obtained by deleting both the $i$th row and the $j$th
column. Let ${\rm {\bf a}}_{.j} $ be the $j$th column and ${\rm {\bf
a}}_{i.} $ be the $i$th row of ${\rm {\bf A}}$. Suppose ${\rm {\bf
A}}_{.j} \left( {{\rm {\bf b}}} \right)$ denotes the matrix obtained from
${\rm {\bf A}}$ by replacing its $j$th column with the column ${\rm {\bf
b}}$, and ${\rm {\bf A}}_{i.} \left( {{\rm {\bf b}}} \right)$ denotes the
matrix obtained from ${\rm {\bf A}}$ by replacing its $i$th row with the
row ${\rm {\bf b}}$.
We  note some properties of column and row determinants of a
quaternion matrix ${\rm {\bf A}} = \left( {a_{ij}} \right)$, where
$i \in I_{n} $, $j \in J_{n} $ and $I_{n} = J_{n} = {\left\{
{1,\ldots ,n} \right\}}$.
\begin{proposition}\label{pr:b_into_brak}  \cite{kyr2}
If $b \in {\mathbb{H}}$, then
\begin{gather*} {\rm{rdet}}_{ i} {\rm {\bf A}}_{i.} \left( {b
\cdot {\rm {\bf a}}_{i.}}  \right) = b \cdot {\rm{rdet}}_{ i} {\rm
{\bf A}},\,\,\,\,{\rm{cdet}} _{{i}}\, {\rm {\bf
A}}_{.i} \left( {{\rm {\bf a}}_{.i}\cdot b} \right) = {\rm{cdet}}
_{{i}}\, {\rm {\bf A}}\cdot  b,\end{gather*} for all $i = 1,\ldots,n $.
\end{proposition}
\begin{proposition} \cite{kyr2}
If for  ${\rm {\bf A}}\in {\rm M}\left( {n,{\mathbb{H}}} \right)$\,
there exists $t \in I_{n} $ such that $a_{tj} = b_{j} + c_{j} $\,
for all $j = 1,\ldots,n$, then
\begin{gather*}
   {\rm{rdet}}_{{i}}\, {\rm {\bf A}} = {\rm{rdet}}_{{i}}\, {\rm {\bf
A}}_{{t{\kern 1pt}.}} \left( {{\rm {\bf b}}} \right) +
{\rm{rdet}}_{{i}}\, {\rm {\bf A}}_{{t{\kern 1pt}.}} \left( {{\rm
{\bf c}}} \right), \,\,\,
  {\rm{cdet}} _{{i}}\, {\rm {\bf A}} = {\rm{cdet}} _{{i}}\, {\rm
{\bf A}}_{{t{\kern 1pt}.}} \left( {{\rm {\bf b}}} \right) +
{\rm{cdet}}_{{i}}\, {\rm {\bf A}}_{{t{\kern 1pt}.}} \left( {{\rm
{\bf c}}} \right),
\end{gather*}
\noindent where ${\rm {\bf b}}=(b_{1},\ldots, b_{n})$, ${\rm {\bf
c}}=(c_{1},\ldots, c_{n})$ and for all ${i = 1,\ldots,n}$.
\end{proposition}
\begin{proposition} \cite{kyr2}
If for ${\rm {\bf A}}\in {\rm M}\left( {n,{\mathbb{H}}} \right)$\,
 there exists $t \in J_{n} $ such that $a_{i\,t} = b_{i} + c_{i}$
for all $i = 1,\ldots,n$, then
\begin{gather*}
  {\rm{rdet}}_{{j}}\, {\rm {\bf A}} = {\rm{rdet}}_{{j}}\, {\rm {\bf
A}}_{.t} \left( {{\rm {\bf b}}} \right) +
{\rm{rdet}}_{{j}}\, {\rm {\bf A}}_{. t} \left(
{{\rm
{\bf c}}} \right),\,\,\,
  {\rm{cdet}} _{j}\, {\rm {\bf A}} = {\rm{cdet}} _{{j}}\, {\rm
{\bf A}}_{.t} \left( {{\rm {\bf b}}} \right) +
{\rm{cdet}} _{j} {\rm {\bf A}}_{.t} \left( {{\rm
{\bf c}}} \right),
\end{gather*}
\noindent where ${\rm {\bf b}}=(b_{1},\ldots, b_{n})^T$, ${\rm
{\bf c}}=(c_{1},\ldots, c_{n})^T$ and for all $j = 1,\ldots,n$.
\end{proposition}
\begin{remark}\label{rem:exp_det}Let $
{\rm{rdet}}_{{i}}\, {\rm {\bf A}} = {\sum\limits_{j = 1}^{n}
{{a_{i{\kern 1pt} j} \cdot R_{i{\kern 1pt} j} } }} $  and $ {\rm{cdet}} _{{j}}\, {\rm {\bf A}} = {{\sum\limits_{i = 1}^{n}
{L_{i{\kern 1pt} j} \cdot a_{i{\kern 1pt} j}} }}$ for all $i,j =
1,\ldots,n$, where by $R_{i{\kern 1pt} j}$ and $L_{i{\kern 1pt} j} $ denote  the right and left $(ij)$th cofactors of ${\rm
{\bf A}}\in {\rm M}\left( {n, {\mathbb{H}}} \right)$, respectively.
It means that  ${\rm{rdet}}_{{i}}\, {\rm
{\bf A}}$ can be expand by right cofactors
  along  the $i$th row and ${\rm{cdet}} _{j} {\rm {\bf A}}$ can be expand by left cofactors
 along  the $j$th column, respectively, for all $i, j = 1,\ldots,n$.
\end{remark}
The main property of the usual determinant is that the determinant of a non-invertible matrix must be  equal zero. But the row and column determinants don't satisfy it, in general. Therefore, these matrix functions  can be consider  as some pre-determinants.
The following theorem has a key value in the theory of the column and row
determinants.
\begin{theorem} \cite{kyr2}\label{theorem:
determinant of hermitian matrix} If ${\rm {\bf A}} = \left( {a_{ij}}
\right) \in {\rm M}\left( {n,{\rm {\mathbb{H}}}} \right)$ is Hermitian,
then ${\rm{rdet}} _{1} {\rm {\bf A}} = \cdots = {\rm{rdet}} _{n} {\rm {\bf
A}} = {\rm{cdet}} _{1} {\rm {\bf A}} = \cdots = {\rm{cdet}} _{n} {\rm {\bf
A}} \in {\rm {\mathbb{R}}}.$
\end{theorem}
Due to Theorem \ref{theorem:
determinant of hermitian matrix}, we can define the
determinant of a  Hermitian matrix ${\rm {\bf A}}\in {\rm M}\left( {n,{\rm
{\mathbb{H}}}} \right)$ by putting,
$\det {\rm {\bf A}}: = {\rm{rdet}}_{{i}}\,
{\rm {\bf A}} = {\rm{cdet}} _{{i}}\, {\rm {\bf A}}, $
 for all $i =1,\ldots,n$.
By using its row and column determinants, the determinant of a quaternion Hermitian matrix   has properties similar to the usual determinant. These properties are completely explored
in
 \cite{kyr2,ky_th}
 and can be summarized in the following
theorems.
\begin{theorem}\label{theorem:row_combin} If the $i$th row of
a Hermitian matrix ${\rm {\bf A}}\in {\rm M}\left( {n,{\rm
{\mathbb{H}}}} \right)$ is replaced with a left linear combination
of its other rows, i.e. ${\rm {\bf a}}_{i.} = c_{1} {\rm {\bf
a}}_{i_{1} .} + \ldots + c_{k}  {\rm {\bf a}}_{i_{k} .}$, where $
c_{l} \in {{\rm {\mathbb{H}}}}$ for all $ l = 1,\ldots,k$ and
$\{i,i_{l}\}\subset I_{n} $, then
\[
 {\rm{rdet}}_{i}\, {\rm {\bf A}}_{i \, .} \left(
{c_{1} {\rm {\bf a}}_{i_{1} .} + \ldots + c_{k} {\rm {\bf
a}}_{i_{k} .}}  \right) = {\rm{cdet}} _{i}\, {\rm {\bf A}}_{i\, .}
\left( {c_{1}
 {\rm {\bf a}}_{i_{1} .} + \ldots + c_{k} {\rm {\bf
a}}_{i_{k} .}}  \right) = 0.
\]
\end{theorem}
\begin{theorem}\label{theorem:colum_combin} If the $j$th column of
 a Hermitian matrix ${\rm {\bf A}}\in
{\rm M}\left( {n,{\rm {\mathbb{H}}}} \right)$   is replaced with a
right linear combination of its other columns, i.e. ${\rm {\bf
a}}_{.j} = {\rm {\bf a}}_{.j_{1}}   c_{1} + \ldots + {\rm {\bf
a}}_{.j_{k}} c_{k} $, where $c_{l} \in{{\rm {\mathbb{H}}}}$ for
all $ l = 1,\ldots,k$ and $\{j,j_{l}\}\subset J_{n}$, then
 \[{\rm{cdet}} _{j}\, {\rm {\bf A}}_{.j}
\left( {{\rm {\bf a}}_{.j_{1}} c_{1} + \ldots + {\rm {\bf
a}}_{.j_{k}}c_{k}} \right) ={\rm{rdet}} _{j} \,{\rm {\bf A}}_{.j}
\left( {{\rm {\bf a}}_{.j_{1}}  c_{1} + \ldots + {\rm {\bf
a}}_{.j_{k}}  c_{k}} \right) = 0.
\]
\end{theorem}
The following theorem about determinantal representation of an
inverse matrix of Hermitian follows immediately from these
properties.
\begin{theorem}
 \cite{ky_th}\label{kyrc10}
 If a Hermitian matrix ${\rm
{\bf A}} \in {\rm M}\left( {n,{\rm {\mathbb{H}}}} \right)$ is such
that $\det {\rm {\bf A}} \ne 0$, then there exist a unique right
inverse  matrix $(R{\rm {\bf A}})^{ - 1}$ and a unique left
inverse matrix $(L{\rm {\bf A}})^{ - 1}$, and $\left( {R{\rm {\bf
A}}} \right)^{ - 1} = \left( {L{\rm {\bf A}}} \right)^{ - 1} =
:{\rm {\bf A}}^{ - 1}$, which possess the following determinantal
representations:
\begin{equation}\label{eq:det_her_inv_R}
  \left( {R{\rm {\bf A}}} \right)^{ - 1} = {\frac{{1}}{{\det {\rm
{\bf A}}}}}
\begin{pmatrix}
  R_{11} & R_{21} & \cdots & R_{n1} \\
  R_{12} & R_{22} & \cdots & R_{n2} \\
  \cdots & \cdots & \cdots & \cdots \\
  R_{1n} & R_{2n} & \cdots & R_{nn}
\end{pmatrix},\end{equation}\begin{equation}\label{eq:det_her_inv_L}
  \left( {L{\rm {\bf A}}} \right)^{ - 1} = {\frac{{1}}{{\det {\rm
{\bf A}}}}}
\begin{pmatrix}
  L_{11} & L_{21} & \cdots  & L_{n1} \\
  L_{12} & L_{22} & \cdots  & L_{n2} \\
  \cdots  & \cdots  & \cdots  & \cdots  \\
  L_{1n} & L_{2n} & \cdots  & L_{nn}
\end{pmatrix},
\end{equation}
where $ \det {\rm {\bf A}} ={\sum\limits_{j = 1}^{n} {{a_{i{\kern
1pt} j} \cdot R_{i{\kern 1pt} j} } }}= {{\sum\limits_{i = 1}^{n}
{L_{i{\kern 1pt} j} \cdot a_{i{\kern 1pt} j}} }}$,

\begin{gather*}
 R_{i j} = {\left\{ {{\begin{array}{*{20}c}
  - {\rm{rdet}}_{{j}}\, {\rm {\bf A}}_{{.{\kern 1pt} j}}^{{i{\kern 1pt} i}} \left( {{\rm
{\bf a}}_{{.{\kern 1pt} {\kern 1pt} i}}}  \right),& {i \ne j},
\hfill \\
 {\rm{rdet}} _{{k}}\, {\rm {\bf A}}^{{i{\kern 1pt} i}},&{i = j},
\hfill \\
\end{array}} } \right.}\,\,\,\,
L_{i j} = {\left\{ {\begin{array}{*{20}c}
 -{\rm{cdet}} _{i}\, {\rm {\bf A}}_{i{\kern 1pt} .}^{j{\kern 1pt}j} \left( {{\rm {\bf a}}_{j{\kern 1pt}. } }\right),& {i \ne
j},\\
 {\rm{cdet}} _{k}\, {\rm {\bf A}}^{j\, j},& {i = j}.
\\
\end{array} }\right.}
\end{gather*}
 The submatrix ${\rm {\bf A}}_{.{\kern 1pt} j}^{i{\kern 1pt} i} \left(
{{\rm {\bf a}}_{.{\kern 1pt} {\kern 1pt} i}}  \right)$ is obtained from
${\rm {\bf A}}$   by  replacing the $j$th column with the $i$th column and then deleting both the $i$th row and column, ${\rm {\bf A}}_{i{\kern 1pt} .}^{jj} \left( {{\rm {\bf
a}}_{j{\kern 1pt} .} } \right)$ is obtained
 by replacing the $i$th row with the $j$th row, and then by
deleting both the $j$th row and  column, respectively. $I_{n} =  {\left\{ {1,\ldots ,n} \right\}}$, $k = \min {\left\{ {I_{n}}
\right.} \setminus {\left. {\{i\}} \right\}}$, for all $i,j =
1,\ldots,n$.
\end{theorem}

\begin{theorem}\label{theorem:equal_lef_righ_her} If ${\rm {\bf A}} \in {\rm M}\left(
{n,{{\rm {\mathbb{H}}}}} \right)$, then $\det {\rm {\bf A}}{\rm
{\bf A}}^{ *} = \det {\rm {\bf A}}^{ * }{\rm {\bf A}}$.
\end{theorem}
\begin{definition}For ${\rm {\bf A}} \in {\rm M}\left(
{n,{{\rm {\mathbb{H}}}}} \right)$, the double determinant of  ${\bf A}$ is defined by putting,
${\rm ddet}\,{\bf A}:= \det {\rm {\bf A}}{\rm
{\bf A}}^{ *} = \det {\rm {\bf A}}^{ * }{\rm {\bf A}}.$

\end{definition}

For  arbitrary ${\rm {\bf A}} \in {\rm M}(n,{{\rm {\mathbb{H}}}})$, we have the following theorem on determinantal representations of its inverse.
\begin{theorem} \label{theorem:deter_inver} The necessary and sufficient condition of invertibility
of  ${\rm {\bf A}} \in {\rm M}(n,{{\rm {\mathbb{H}}}})$ is
${\rm{ddet}} {\rm {\bf A}} \ne 0$. Then there exists ${\rm {\bf
A}}^{ - 1} = \left( {L{\rm {\bf A}}} \right)^{ - 1} = \left(
{R{\rm {\bf A}}} \right)^{ - 1}$, where
\begin{equation*}
 \left( {L{\rm {\bf A}}} \right)^{ - 1}
=\left( {{\rm {\bf A}}^{ *}{\rm {\bf A}} } \right)^{ - 1}{\rm {\bf
A}}^{ *} ={\frac{{1}}{{{\rm{ddet}}{ \rm{\bf A}} }}}
\begin{pmatrix}
  {\mathbb{L}} _{11} & {\mathbb{L}} _{21}& \ldots & {\mathbb{L}} _{n1} \\
  {\mathbb{L}} _{12} & {\mathbb{L}} _{22} & \ldots & {\mathbb{L}} _{n2} \\
  \ldots & \ldots & \ldots & \ldots \\
 {\mathbb{L}} _{1n} & {\mathbb{L}} _{2n} & \ldots & {\mathbb{L}} _{nn}
\end{pmatrix}
\end{equation*}
\begin{equation*} \left( {R{\rm {\bf A}}} \right)^{ - 1} = {\rm {\bf
A}}^{ *} \left( {{\rm {\bf A}}{\rm {\bf A}}^{ *} } \right)^{ - 1}
= {\frac{{1}}{{{\rm{ddet}}{ \rm{\bf A}}^{ *} }}}
\begin{pmatrix}
 {\mathbb{R}} _{11} & {\mathbb{R}} _{21} &\ldots & {\mathbb{R}} _{n1} \\
 {\mathbb{R}} _{12} & {\mathbb{R}} _{22} &\ldots & {\mathbb{R}} _{n2}  \\
 \ldots  & \ldots & \ldots & \ldots \\
 {\mathbb{R}} _{1n} & {\mathbb{R}} _{2n} &\ldots & {\mathbb{R}} _{nn}
\end{pmatrix}
\end{equation*}
and \[{\mathbb{L}} _{ij} = {\rm{cdet}} _{j} ({\rm {\bf
A}}^{\ast}{\rm {\bf A}})_{.j} \left( {{\rm {\bf a}}_{.{\kern 1pt}
i}^{ *} } \right), \,\,\,{\mathbb{R}} _{\,{\kern 1pt} ij} =
{\rm{rdet}}_{i} ({\rm {\bf A}}{\rm {\bf A}}^{\ast})_{i.} \left(
{{\rm {\bf a}}_{j.}^{ *} }  \right),\] for all $i,j =
{1,\ldots,n}.$
\end{theorem}
Moreover, the following  criterion of invertibility of a quaternion matrix
can be obtained.
\begin{theorem}\label{theorem:inver_equiv} If ${\rm {\bf A}}\in {\rm M}\left( {n,{\rm
{\mathbb{H}}}} \right)$, then
 the following statements are equivalent.
\begin{itemize}
  \item[i)] ${\rm {\bf A}}$ is invertible, i.e. ${\rm {\bf A}}\in
GL\left( {n,{\mathbb{H}}} \right);$
  \item[ii)]  rows of ${\rm {\bf A}}$ are left-linearly independent;
  \item[iii)]  columns of ${\rm {\bf A}}$ are right-linearly
independent;
  \item[iv)] ${\rm ddet}\, {\rm {\bf A}}\ne 0$.
\end{itemize}
\end{theorem}
\subsection{Some provisions of quaternion eigenvalues}
Due to real-scalar multiplying on the right, quaternion column-vectors  form a right  vector $\mathbb{R}$-space, and, by real-scalar multiplying on the left, quaternion row-vectors  form a left  vector $\mathbb{R}$-space denoted by $\mathcal{H}_{r}$ and $\mathcal{H}_{l}$, respectively. It can be shown that $\mathcal{H}_{r}$ and $\mathcal{H}_{l}$ possess corresponding $\mathbb{H}$-valued inner products by putting $\langle \mathbf{x},\mathbf{y}\rangle_{r}=\overline{y}_1x_{1}+\cdots+\overline{y}_{n}x_{n}$ for $\mathbf{x}=\left(x_{i}\right)_{i=1}^{n}, \mathbf{y}=\left(y_{i}\right)_{i=1}^{n}\in \mathcal{H}_{r}$, and $\langle \mathbf{x},\mathbf{y}\rangle_{l}=x_{1}\overline{y}_1+\cdots+x_{n}\overline{y}_{n}$ for $\mathbf{x}, \mathbf{y}\in \mathcal{H}_{l}$ that satisfy the  inner product relations, namely, conjugate symmetry, linearity, and positive-definiteness but with specialties
 \begin{equation*}\begin{gathered}\langle \mathbf{x}\alpha+\mathbf{y}\beta,\mathbf{z} \rangle=\langle \mathbf{x},\mathbf{z}\rangle\alpha+\langle \mathbf{y},\mathbf{z}\rangle\beta \,\,\mbox{ when}\,\, \mathbf{x}, \mathbf{y}, \mathbf{z} \in \mathcal{H}_{r}\\
          \langle \alpha\mathbf{x}+\beta\mathbf{y},\mathbf{z} \rangle=\alpha\langle \mathbf{x},\mathbf{z}\rangle+\beta\langle \mathbf{y},\mathbf{z}\rangle \,\, \mbox{ when}\,\, \mathbf{x}, \mathbf{y}, \mathbf{z} \in \mathcal{H}_{l},\end{gathered}\end{equation*}
 for any $\alpha, \beta \in {\mathbb{H}}$.
 A set of vectors from  $\mathcal{H}_{r}$ and $\mathcal{H}_{l}$ can be orthonormalize in particular by the Gram–Schmidt process with corresponding projection operators
  \begin{gather*}
 {\rm proj}_{\bf u}({\bf v}):= {\bf u}\frac{\langle \mathbf{u},\mathbf{v}\rangle_{r}}{\langle \mathbf{u},\mathbf{u}\rangle_{r}},\\
 {\rm proj}_{\bf u}({\bf v}):= \frac{\langle \mathbf{u},\mathbf{v}\rangle_{l}}{\langle \mathbf{u},\mathbf{u}\rangle_{l}}{\bf u}
 \end{gather*}
 for $\mathcal{H}_{r}$ and $\mathcal{H}_{l}$, respectively.
 Due to the above, the following definition makes sense.
 \begin{definition}
Suppose ${\rm {\bf U}} \in {\rm M}\left( {n,{\rm {\mathbb{H}}}}
\right)$ and ${\rm {\bf U}}^{ *} {\rm {\bf U}} = {\rm {\bf U}}{\rm
{\bf U}}^{ *}  = {\rm {\bf I}}$, then the matrix ${\rm {\bf U}}$
is called unitary.
\end{definition}
Clear, that columns of ${\bf U}$ form a system of normalized vectors in  $\mathcal{H}_{r}$, rows of ${\bf U}^{*}$ is a system of normalized vectors in  $\mathcal{H}_{l}$.

 Due to the noncommutativity of quaternions, there are two
types of eigenvalues.
A quaternion $\lambda$ is said to be a left eigenvalue
of ${\rm {\bf A}} \in {\rm M}\left( {n,{\rm {\mathbb{H}}}}
\right)$ if
$
{\rm
{\bf A}} \cdot {\rm {\bf x}} = \lambda \cdot {\rm {\bf x}}$,
and a right eigenvalue if
 $
{\rm {\bf A}} \cdot {\rm {\bf x}} = {\rm
{\bf x}} \cdot \lambda$
for some nonzero quaternion column-vector
${\rm {\bf x}}$. Then, the  set  $\{\lambda\in {\mathbb{H}}|  {\bf A}  {\rm {\bf x}} = \lambda  {\rm {\bf x}},\,  {\bf x}\neq {\bf 0}  \in {\mathbb{H}}^{n}\} $ is  called  the
left  spectrum  of  ${\bf A}$,  denoted  by  $\sigma_{l}({\bf A})$.  The  right  spectrum  is  similarly  defined
  by putting, $\sigma_{r}({\bf A}):=\{\lambda\in {\mathbb{H}}|  {\bf A} {\rm {\bf x}} ={\rm {\bf x}}\lambda,\,  {\bf x}\neq {\bf 0}  \in {\mathbb{H}}^{n}\} $.

The theory on the left eigenvalues of quaternion matrices has been
investigated in particular in \cite{hu, so, wo}. The theory on the
right eigenvalues of quaternion matrices is more developed \cite{br,ma,ba,dra,zh,far}.
 We consider this is a natural consequence of the fact that quaternion column vectors form a right  vector space for which left eigenvalues  seem to be "exotic" because of their multiplying from the left.

 We present the some known results from the theory of right eigenvalues.
 It's well known that if   $\lambda$ is  a  nonreal  eigenvalue  of   ${\bf A}$,  so  is  any  element  in the
equivalence  class  containing   $[\lambda]$, i.e. $[\lambda]=\{x|x= u^{-1}\lambda u,\, u\in{\mathbb{H}},\, \|u\|=1\}$.
 \begin{theorem}\cite{br} Any   quaternion  matrix  ${\rm {\bf A}} \in {\rm M}\left( {n,{\rm {\mathbb{H}}}}
\right)$  has  exactly  $n$    eigenvalues  which  are  complex  numbers  with
nonnegative  imaginary    parts.
\end{theorem}
Those  eigenvalues $h_{1} + k_{1}{\bf  i},\ldots, h_{n} + k_{n}{\bf  i}$, where $k_{t}\geq0$ and
$h_{t}, k_{t}\in {\mathbb{R}}$  for all $t = 1,\ldots, n$, are  said  to  be  the  standard  eigenvalues  of  ${\bf A}$.
\begin{theorem}\cite{br} Let $ {\bf A} \in {\rm M}\left( n, {\mathbb{H}}
\right)$. Then there exists a unitary matrix ${\bf U}$ such that ${\bf U}^{*}{\bf AU}$ is an upper
triangular matrix with diagonal entries $h_{1} + k_{1}{\bf  i},\ldots, h_{n} + k_{n}{\bf  i}$ which are the standard eigenvalues of  ${\bf A}$.
\end{theorem}
\begin{corollary}\cite{zh}\label{th:diag_spec}
Let ${\rm {\bf A}} \in {\rm M}\left( {n,{\rm {\mathbb{H}}}}
\right)$ with the standard eigenvalues $h_{1} + k_{1}{\bf  i},\ldots, h_{n} + k_{n}{\bf i}$.
Then
$
\sigma_{r}=[h_{1} + k_{1}{\bf  i}]\cup\cdots\cup[h_{n} + k_{n}{\bf  i}].
$
\end{corollary}

 \begin{corollary}\cite{zh}\label{th:decomp_herm}
Let ${\rm {\bf A}} \in {\rm M}\left( {n,{\rm {\mathbb{H}}}}
\right)$ be given. Then, ${\rm {\bf A}}$ is Hermitian  if and only
if there are a unitary matrix ${\rm {\bf U}} \in {\rm M}\left(
{n,{\rm {\mathbb{H}}}} \right)$ and a real diagonal matrix ${\rm
{\bf D}} = {\rm diag}\left( {\lambda _{{\kern 1pt} 1} ,\lambda
_{{\kern 1pt} 2} ,\ldots ,\lambda _{{\kern 1pt} n}}  \right)$ such
that ${\rm {\bf A}} = {\rm {\bf U}}{\rm {\bf D}}{\rm {\bf U}}^{
*}$, where $\lambda _{ 1},...,\lambda _{ n} $ are right
eigenvalues of ${\rm {\bf A}}$.
\end{corollary}

The right  and left  eigenvalues are in general unrelated \cite{fa}, but it is not  for Hermitian matrices.
Suppose ${\rm {\bf A}} \in {\rm M}\left( {n, {\mathbb{H}}}\right)$
is Hermitian and $\lambda \in {\rm {\mathbb {R}}}$ is its right
eigenvalue, then ${\rm {\bf A}} \cdot {\rm {\bf x}} = {\rm {\bf
x}} \cdot \lambda = \lambda \cdot {\rm {\bf x}}$. This means that
all right eigenvalues of a Hermitian matrix are its left
eigenvalues as well. For real left eigenvalues, $\lambda \in {\rm
{\mathbb {R}}}$, the matrix $\lambda {\rm {\bf I}} - {\rm {\bf
A}}$ is Hermitian.

\begin{definition}
If $t \in {\rm {\mathbb {R}}}$, then for a Hermitian matrix ${\rm
{\bf A}}$ the polynomial $p_{{\rm {\bf A}}}\left( {t} \right) =
\det \left( {t{\rm {\bf I}} - {\rm {\bf A}}} \right)$ is said to
be the characteristic polynomial of ${\rm {\bf A}}$.
\end{definition}

\begin{lemma}\cite{ky_th}\label{lem:char_her}
If ${\rm {\bf A}} \in {\rm M}\left( {n,{\rm {\mathbb{H}}}}
\right)$ is Hermitian, then $p_{{\rm {\bf A}}}\left( {t} \right) =
t^{n} - d_{1} t^{n - 1} + d_{2} t^{n - 2} - \ldots + \left( { - 1}
\right)^{n}d_{n}$, where $d_{k} $ is the sum of principle minors
of ${\rm {\bf A}}$ of order $k$, $1 \le k < n$, and $d_{n}=\det
{\rm {\bf A}}$.
\end{lemma}

The roots of the
characteristic polynomial of a Hermitian matrix are its real left
eigenvalues, which are its right eigenvalues as well.

\subsection{ Determinantal representations of the Moore-Penrose and weighted  Moore-Penrose inverses over the quaternion skew field}

Within the framework of the theory of   column-row determinants, we have the following theorem on determinantal representations of the quaternion Moore-Penrose inverse.
\begin{theorem} \cite{ky_math_sci}\label{theor:det_repr_MP}
If ${\rm {\bf A}} \in {\rm {\mathbb{H}}}_{r}^{m\times n} $, then
the Moore-Penrose inverse  ${\rm {\bf A}}^{ \dag} = \left( {a_{ij}^{
\dag} } \right) \in {\rm {\mathbb{H}}}_{}^{n\times m} $ possess the
following determinantal representations:
\begin{enumerate}
  \item[(i)] If $r< {\rm min} \{m, n\}$, then
  \begin{equation}
\label{eq:det_repr_A*A}
 a_{ij}^{ \dag}  = {\frac{{{\sum\limits_{\beta
\in J_{r,\,n} {\left\{ {i} \right\}}} {{\rm{cdet}} _{i} \left(
{\left( {{\rm {\bf A}}^{ *} {\rm {\bf A}}} \right)_{\,. \,i}
\left( {{\rm {\bf a}}_{.j}^{ *} }  \right)} \right){\kern 1pt}
{\kern 1pt} _{\beta} ^{\beta} } } }}{{{\sum\limits_{\beta \in
J_{r,\,\,n}} {{{\rm{det}} {\left( {{\rm {\bf A}}^{ *} {\rm {\bf A}}}
\right){\kern 1pt} _{\beta} ^{\beta} }  }}} }}},
\end{equation}
or
  \begin{equation}
\label{eq:det_repr_AA*} a_{ij}^{ \dag}  =
{\frac{{{\sum\limits_{\alpha \in I_{r,m} {\left\{ {j} \right\}}}
{{\rm{rdet}} _{j} \left( {({\rm {\bf A}}{\rm {\bf A}}^{ *}
)_{j\,.\,} ({\rm {\bf a}}_{i.\,}^{ *} )} \right)\,_{\alpha}
^{\alpha} } }}}{{{\sum\limits_{\alpha \in I_{r,\,m}}  {{{\rm{det}}
{\left( {{\rm {\bf A}}{\rm {\bf A}}^{ *} } \right){\kern 1pt}
_{\alpha} ^{\alpha} } }}} }}}.
\end{equation}
 \item [(ii)]   If $r= n$, then
   \begin{equation}
\label{eq:det_repr_f_A*A}
a_{ij}^{ \dag}  =  {\frac{{{\rm cdet} _{i} ({\rm
{\bf A}}^{ *} {\rm {\bf A}})_{.\,i} \left( {{\rm {\bf a}}_{.j}^{ *} }  \right)}}{{ {\rm det} ({\rm
{\bf A}}^{ *} {\rm {\bf A}})}}}
  \end{equation}
  or (\ref{eq:det_repr_AA*}) when $n< m$.

  \item [(iii)]
  If $r= m$, then
   \begin{equation}
\label{eq:det_repr_f_AA*}
 a_{ij}^{\dag}  = {\frac{{{\rm rdet} _{j} ({\rm
{\bf A}} {\rm {\bf A}}^{ *})_{j\,.} \left( {{\rm {\bf a}}_{i.}^{ *} }  \right)}}{{ {\rm det} ({\rm
{\bf A}} {\rm {\bf A}}^{ *})}}}
  \end{equation}
  or (\ref{eq:det_repr_A*A}) when $m<n $.
\end{enumerate}
\end{theorem}

Even though  the  eigenvalues of ${\bf A}^{\sharp}{\bf A}$ and ${\bf A}{\bf A}^{\sharp}$ are real and nonnegative, they are not Hermitian in general. Therefor, the following two cases are considered, when ${\bf A}^{\sharp}{\bf A}$ and ${\bf A}{\bf A}^{\sharp}$ both or one of them are Hermitian, and when they are non-Hermitian. Denote the $(ij)$th entry of ${\rm {\bf A}}_{M,N}^{\dag}$ by $a_{ij}^{\ddag}$ for all $i=1,\ldots,n$ and $j=1,\ldots,m$.
\begin{theorem}\cite{kyr_wmpi}\label{th:det_rep_A_mn}
Let ${\rm {\bf A}} \in {\rm {\mathbb{H}}}_{r}^{m\times n} $. If  ${\bf A}^{\sharp}{\rm {\bf A}}$ or ${\bf A}{\rm {\bf A}}^{\sharp}$ are Hermitian, then
the weighted Moore-Penrose inverse  ${\rm {\bf A}}_{M,N}^{\dag} = \left( {{a}_{ij}^{\ddag} } \right) \in
{\rm {\mathbb{H}}}^{n\times m} $ possess the following determinantal representations, respectively,
\begin{enumerate}
  \item[(i)] If $r< {\rm min} \{m, n\}$, then
\begin{equation}
\label{kyr5} {a}_{ij}^{\ddag}  = {\frac{{{\sum\limits_{\beta \in
J_{r,\,n} {\left\{ {i} \right\}}} {{\rm{cdet}} _{i} \left( {\left(
{{\rm {\bf A}}^{\sharp} {\rm {\bf A}}} \right)_{\,. \,i} \left( {{\rm
{\bf a}}_{.j}^{\sharp} }  \right)} \right){\kern 1pt} {\kern 1pt}
_{\beta} ^{\beta} } } }}{{{\sum\limits_{\beta \in J_{r,\,\,n}}
{{\left| {\left( {{\rm {\bf A}}^{\sharp} {\rm {\bf A}}} \right){\kern
1pt}  _{\beta} ^{\beta} }  \right|}}} }}},
\end{equation}
or
\begin{equation}
\label{kyr6} {a}_{ij}^{\ddag}  = {\frac{{{\sum\limits_{\alpha \in
I_{r,m} {\left\{ {j} \right\}}} {{\rm{rdet}} _{j} \left( {({\rm
{\bf A}}{\rm {\bf A}}^{\sharp} )_{j\,.\,} ({\rm {\bf a}}_{i.\,}^{\sharp}
)} \right)\,_{\alpha} ^{\alpha} } }}}{{{\sum\limits_{\alpha \in
I_{r,\,m}}  {{\left| {\left( {{\rm {\bf A}}{\rm {\bf A}}^{\sharp} }
\right){\kern 1pt}  _{\alpha} ^{\alpha} } \right|}}} }}}.
\end{equation}
 \item [(ii)]   If $\rank{\rm {\bf A}} = n<m$, then
\begin{equation}\label{eq:a_+a1a_f}
 {{a}_{ij}^{\ddag} }   = {\frac{{{\rm{cdet}} _{i} ({{\rm {\bf A}}^{\sharp} {\rm {\bf
A}}})_{.\,i} \left( {{\rm {\bf a}}_{.\,j}^{\sharp} }  \right)} }{{ {\rm det} ({{\rm {\bf A}}^{\sharp} {\rm {\bf
A}}})}}},
\end{equation}
or the determinantal representation (\ref{kyr5}) can be applicable  as well.
 \item [(iii)]  If $\rank{\rm {\bf A}} = m<n$,  then
 \begin{equation}
\label{eq:a_+aa1_f}  {{a}_{ij}^{\ddag} }   =    {\frac{{{\rm{rdet}} _{j} ({\rm {\bf A}} {\rm {\bf A}}^{\sharp})_{j.} \left( {{\rm {\bf a}}_{i.}^{\sharp} }  \right)}}{{{\rm det} {({\rm {\bf A}} {\rm {\bf A}}^{\sharp}) }}}}.
\end{equation}
or the determinantal representation (\ref{kyr6}) can be applicable  as well.
\end{enumerate}
\end{theorem}
Denote ${\bf M}^{\frac{1}{2}}=\left(m_{ij}^{(\frac{1}{2})}\right)$, ${\bf N}^{-\frac{1}{2}}=\left(n_{ij}^{(-\frac{1}{2})}\right)$, and $\widetilde{{\bf A}}:=
{\bf M}^{\frac{1}{2}}{\bf A}{\bf N}^{-\frac{1}{2}}
 =\left(\widetilde{a}_{ij}\right)\in {\mathbb{H}}^{m\times n} $, then ${\bf N}^{-\frac{1}{2}}{\bf A}^{*}{\bf M}^{\frac{1}{2}}=\widetilde{{\bf A}}^{*}=\left(\widetilde{a}^{*}_{ij}\right)$, $\left({\bf M}^{\frac{1}{2}}{\bf A}{\bf N}^{-\frac{1}{2}}\right)^{\dag}=\widetilde{{\bf A}}^{\dag}=\left(\widetilde{a}^{\dag}_{ij}\right)$.
 \begin{theorem}\cite{kyr_wmpi}\label{th:det_rep_A_mn_non}
Let ${\rm {\bf A}} \in {\rm {\mathbb{H}}}_{r}^{m\times n} $.
\begin{enumerate}
  \item[(i)] If ${\bf A}^{\sharp}{\rm {\bf A}}$ is non-Hermitian, then
the weighted Moore-Penrose inverse  ${\rm {\bf A}}_{M,N}^{\dag} = \left( {{a}_{ij}^{\ddag} } \right) \in
{\rm {\mathbb{H}}}^{n\times m} $ possess the  determinantal representations
\begin{enumerate}
  \item if $r<n$
  \begin{equation}\label{eq:det_rep_non_c1}
{{a}_{ij}^{\ddag} }=
{\frac{\sum\limits_{k}n^{(-\frac{1}{2})}_{ik}{{\sum\limits_{\beta \in
J_{r,\,n} {\left\{ {i} \right\}}} {{\rm{cdet}} _{k} \left( {\left(
\widetilde{{\bf A}}^{*}\widetilde{{\bf A}}
 \right)_{\,. \,k}  \left({\widehat{\bf a}}_{.j}\right)  } \right)
_{\beta} ^{\beta} } } }     }{{{\sum\limits_{\beta \in J_{r,\,\,n}}
{{\left| {\left(
\widetilde{{\bf A}}^{*}\widetilde{{\bf A}}
 \right)  _{\beta} ^{\beta} }  \right|}}} }}},\end{equation}
 where $\widehat{{\bf a}}_{.j}$ is the $j$th column of ${\bf N}^{-\frac{1}{2}}{\bf A}^{*}{\bf M}$;
  \item if $r=n$
  \begin{equation} \label{eq:det_rep_non_c_ful1}
{{a}_{ij}^{\ddag} }=\frac{{\rm{cdet}} _{i}({\bf A}^{*}{\bf M}{\bf A})_{.i}\,(\widehat{\bf{a}}_{.j})}{\det({\bf A}^{*}{\bf M}{\bf A})},
\end{equation}
where $\widehat{\bf{a}}_{.j}$ is the $j$th column of ${\bf A}^{*}{\bf M}$ for all $j=1,\ldots,m$.

\end{enumerate}
  \item[(ii)] If ${\bf A}{\rm {\bf A}}^{\sharp}$ is non-Hermitian , then
  ${\rm {\bf A}}_{M,N}^{\dag} = \left( {{a}_{ij}^{\ddag} } \right) $ possess the  determinantal representation
   \begin{enumerate}
     \item if $r<m$, \begin{equation}\label{eq:det_rep_non_r1}
{{a}_{ij}^{\ddag} }=
 {\frac{\sum\limits_{l}{{\sum\limits_{\alpha \in
I_{r,\,m} {\left\{ {l} \right\}}} {{\rm{rdet}} _{l} \left( {\left(
\widetilde{{\bf A}}\widetilde{{\bf A}}^{*}
 \right)_{l. \,}  (\widehat{{\bf a}}_{i.})  } \right)
_{\alpha} ^{\alpha} } } }\cdot m^{(\frac{1}{2})}_{lj}}{{{\sum\limits_{\alpha \in I_{r,\,\,m}}
{{\left| {\left( \widetilde{{\bf A}}\widetilde{{\bf A}}^{*} \right) _{\alpha} ^{\alpha} }  \right|}}} }}},
\end{equation}
where $\widehat{{\bf a}}_{i.}$ is the $i$th row of ${\bf N}^{-1}{\bf A}^{*}{\bf M}^{\frac{1}{2}}$;
     \item if $r=m$, \begin{equation}\label{eq:det_rep_non_r_ful1}
{{a}_{ij}^{\ddag} }=\frac{{\rm{rdet}}_{j}({\bf A}{\bf N}^{-1}{\bf A}^{*})_{j.}(\widehat{\bf{a}}_{i.})}{\det({\bf A}{\bf N}^{-1}{\bf A}^{*})}.
\end{equation}
where $\widehat{\bf{a}}_{i.}$ is the $i$th row of ${\bf N}^{-1}{\bf A}^{*}$ for all $i=1,\ldots,n$.
 \end{enumerate}
  \end{enumerate}
 \end{theorem}

\section{Cramer's Rule for Two-sided Restricted\\  Quaternionic Matrix Equation}
\begin{definition}
For an arbitrary matrix over the quaternion skew field, ${\bf A}\in  {\mathbb{H}}^{m\times n}$, we denote by
\begin{itemize}
  \item  $\mathcal{R}_{r}({\rm {\bf A}})=\{ {\bf y}\in {\mathbb{H}}^{m\times 1} : {\bf y} = {\bf A}{\bf x},\,  {\bf x} \in {\mathbb{H}}^{n\times 1}\},$  the column right space of ${\bf A}$,
  \item  $\mathcal{N}_{r}({\rm {\bf A}})=\{ {\bf x}\in {\mathbb{H}}^{n\times 1} : \,\, {\bf A}{\bf x}=0\}$,  the right null space of  ${\bf A}$,
  \item $\mathcal{R}_{l}({\rm {\bf A}})=\{ {\bf y}\in {\mathbb{H}}^{1\times n} : \,\,{\bf y} = {\bf x}{\bf A},\,\,  {\bf x} \in {\mathbb{H}}^{1\times m}\}$, the row left space of ${\bf A}$,
  \item  $\mathcal{N}_{l}({\rm {\bf A}})=\{ {\bf x}\in {\mathbb{H}}^{1\times m} : \,\, {\bf x}{\bf A}=0\}$,  the left null space of  ${\bf A}$.
\end{itemize}
\end{definition}

It is easy to see, if ${\bf A}\in  {\mathbb{H}}_{n}^{n\times n}$, then
$
\mathcal{R}_{r}\oplus\mathcal{N}_{r}={\mathbb{H}}^{n\times 1},$ and $\mathcal{R}_{l}\oplus\mathcal{N}_{l}={\mathbb{H}}^{1\times n}.
$
Suppose that ${\bf A}\in {\mathbb{H}}^{m\times n}$, ${\bf B}\in {\mathbb{H}}^{p\times q}$. Denote
 \begin{multline*}
 \mathcal{R}_{r}( {\bf A},{\bf B}):=\mathcal{N}_{r}({\bf  Y})=\{{\bf  Y}={\bf A}{\bf X}{\bf B}:{\bf X}^{n\times q}\},\\
\mathcal{N}_{r}( {\bf A},{\bf B}):=\mathcal{R}_{r}({\bf  X})=\{{\bf X}^{n\times p}:{\bf A}{\bf X}{\bf B}={\bf 0}\},\\
\mathcal{R}_{l}( {\bf A},{\bf A}:)=\mathcal{R}_{l}({\bf  Y})=\{{\bf  Y}={\bf A}{\bf X}{\bf B}:{\bf X}^{n\times q}\},\\
 \mathcal{N}_{l}( {\bf A},{\bf B}):=\mathcal{N}_{l}({\bf  X})=\{{\bf X}^{n\times p}:{\bf A}{\bf X}{\bf B}={\bf 0}\}.
\end{multline*}

\begin{lemma}\label{lem:sol_res}\cite{song6} Suppose that ${\bf A}\in {\mathbb{H}}^{m\times n}_{{r}_1}$, ${\bf B}\in {\mathbb{H}}^{p\times q}_{{r}_2}$,  ${\bf M}$, ${\bf N}$, ${\bf P}$, and ${\bf Q}$ are Hermitian positive definite matrices of order $m$, $n$, $p$, and $q$, respectively. Denote ${\bf A}^{\sharp}={\bf N}^{-1}{\bf A}^{*}{\bf M}$ and ${\bf B}^{\sharp}={\bf Q}^{-1}{\bf B}^{*}{\bf P}$. If ${\bf D}\subset  \mathcal{R}_{r}\left({\bf A}{\bf A}^{\sharp}, {\bf B}^{\sharp}{\bf B}\right)$ and ${\bf D}\subset  \mathcal{R}_{l}\left({\bf A}^{\sharp}{\bf A}, {\bf B}{\bf B}^{\sharp}\right)$,
\begin{gather}\label{eq1}
{\bf A}{\bf X}{\bf B}={\bf D},\\\label{res1}
\mathcal{R}_{r}({\bf X})\subset {\bf N}^{-1}\mathcal{R}_{r}({\bf A}^{*}),\,\mathcal{N}_{r}({\bf X})\supset {\bf P}^{-1}\mathcal{N}_{r}({\bf B}^{*}),\\\label{res2}\mathcal{R}_{l}({\bf X})\subset \mathcal{R}_{l}({\bf A}^{*}){\bf M},\,\mathcal{N}_{l}({\bf X})\supset \mathcal{N}_{l}({\bf B}^{*}){\bf Q}
\end{gather}
then the unique solution of (\ref{eq1}) with the restrictions (\ref{res1})-(\ref{res2}) is
\begin{equation}\label{eq:sol_res}
{\bf X}= {\bf A}^{\dag}_{M,N}{\bf D}{\bf B}^{\dag}_{P,Q}.\end{equation}
\end{lemma}
In this chapter, we get determinantal representations of (\ref{eq:sol_res}) that are intrinsically analogs of the classical Cramer's rule. We will consider several cases depending on whether the matrices ${\bf A}^{\sharp}{\bf A}$ and ${\bf B}{\bf B}^{\sharp}$ are Hermitian or not.

\subsection{The Case of Both Hermitian Matrices ${\bf A}^{\sharp}{\bf A}$ and ${\bf B}{\bf B}^{\sharp}$.}
Denote
${\rm {\bf \widetilde{D}}}= {\rm {\bf
A}}^{\sharp}{\rm {\bf D}}{\rm {\bf B}}^{\sharp}$.
\begin{theorem}\label{th:cr_rul_AXB} Let ${\bf A}^{\sharp}{\bf A}$ and ${\bf B}{\bf B}^{\sharp}$ be Hermitian.  Then the  solution  (\ref{eq:sol_res})  possess the following determinantal representations.
\begin{enumerate}
\item[(i)] If ${\rm rank}\,{\rm {\bf A}} = r_{1} <  n$ and
 ${\rm rank}\,{\rm {\bf B}} = r_{2} < p$, then
\begin{equation}\label{eq:d^B}
x_{ij} = {\frac{{{\sum\limits_{\beta \in J_{{{r}_1},\,n} {\left\{
{i} \right\}}} {{\rm{cdet}} _{i} \left( {\left( {{\rm {\bf A}}^{\sharp} {\rm {\bf A}}} \right)_{.\,i} \left( {{{\rm {\bf
d}}}\,_{.\,j}^{{\rm {\bf B}}}} \right)} \right) _{\beta} ^{\beta}
} } }}{{{\sum\limits_{\beta \in J_{{{r}_1},n}} {{\left| {\left(
{{\rm {\bf A}}^{\sharp} {\rm {\bf A}}} \right)_{\beta} ^{\beta} }
\right|}} \sum\limits_{\alpha \in I_{{{r}_2},p}}{{\left| {\left(
{{\rm {\bf B}}{\rm {\bf B}}^{\sharp} } \right) _{\alpha} ^{\alpha} }
\right|}}} }}},
\end{equation} or

\begin{equation}\label{eq:d^A}
 x_{ij}={\frac{{{\sum\limits_{\alpha
\in I_{r_{2},p} {\left\{ {j} \right\}}} {{\rm{rdet}} _{j} \left(
{\left( {{\rm {\bf B}}{\rm {\bf B}}^{\sharp} } \right)_{\,j\,.} \left(
{{{\rm {\bf d}}}\,_{i\,.}^{{\rm {\bf A}}}} \right)}
\right)\,_{\alpha} ^{\alpha} } }}}{{{\sum\limits_{\beta \in
J_{r_{1},n}} {{\left| {\left( {{\rm {\bf A}}^{\sharp} {\rm {\bf A}}}
\right) _{\beta} ^{\beta} } \right|}}\sum\limits_{\alpha \in
I_{r_{2},p}} {{\left| {\left( {{\rm {\bf B}}{\rm {\bf B}}^{\sharp} }
\right) _{\alpha} ^{\alpha} } \right|}}} }}},
\end{equation}

where
\begin{equation} \label{eq:def_d^B_m}
   {{{\rm {\bf d}}}_{.\,j}^{{\rm {\bf B}}}}=\left(
\sum\limits_{\alpha \in I_{{{r}_2},p} {\left\{ {j} \right\}}}
{{\rm{rdet}} _{j} \left( {\left( {{\rm {\bf B}}{\rm {\bf B}}^{\sharp}
} \right)_{j.} \left( {\tilde{{\rm {\bf d}}}_{k.}} \right)}
\right)_{\alpha} ^{\alpha}}
\right)\in{\rm {\mathbb{H}}}^{n \times
1}\end{equation}
\begin{equation} \label{eq:d_A<n}  {{{\rm {\bf d}}}_{i\,.}^{{\rm {\bf A}}}}=\left(
\sum\limits_{\beta \in J_{r_{1},n} {\left\{ {i} \right\}}}
{{\rm{cdet}} _{i} \left( {\left( {{\rm {\bf A}}^{\sharp}{\rm {\bf A}}
} \right)_{.i} \left( {\tilde{{\rm {\bf d}}}_{.l}} \right)}
\right)_{\beta} ^{\beta}} \right)\in{\rm {\mathbb{H}}}^{1 \times
p}
\end{equation}
are the column-vector and the row-vector, respectively.  ${\tilde{{\rm {\bf
 d}}}_{k.}}$ and
${\tilde{{\rm {\bf d}}}_{.l}}$ are the $k$th row   and the $l$th
column  of ${\rm {\bf \widetilde{D}}}$ for all $k =1,...,n $, $l =1,...,p $.
\item[(ii)] If ${\rm rank}\,{\rm {\bf A}} = n$ and ${\rm rank}\,{\rm {\bf B}} = p$, then
\begin{equation}
\label{eq:AXB_cdetA*A} x_{i\,j} = \frac{{{\rm cdet} _{i} (
{\bf A}^{\sharp}  {\bf A})_{.\,i\,} \left( {{\rm {\bf
d}}_{.j}^{{\rm {\bf B}}}} \right)}}{{\rm det} ({\bf A}^{\sharp}  {\bf A})\cdot {\rm det} (
{\bf B} {\bf B}^{\sharp} )},
\end{equation}
or
\begin{equation}
\label{eq:AXB_rdetBB*} x_{i\,j} = \frac{{\rm rdet} _{j} (
{\bf B} {\bf B}^{\sharp} )_{j.\,} \left(  {\bf
d}_{i\,.}^{ {\bf A}} \right)}{{\rm det} ({\bf A}^{\sharp}  {\bf A})\cdot {\rm det} (
{\bf B} {\bf B}^{\sharp} )},
\end{equation}
 \noindent where   \begin{equation}\label{eq:d_B_p} {\bf d}_{.j}^ {\bf B} : =
\left( {\rm
rdet} _{j} ( {\bf B} {\bf B}^{\sharp} )_{j.\,} \left(
{\tilde {\bf d}}_{k\,.} \right) \right)\in{\mathbb{H}}^{n \times
1},\end{equation}
\begin{equation}\label{eq:def_d^A}{\bf d}_{i\,.}^ {\bf A} : = \left(
{\rm cdet}
_{i} ( {\bf A}^{\sharp}  {\bf A})_{.\,i} \left(
{\tilde {\bf d}}_{.l} \right) \right)\in{\rm {\mathbb{H}}}^{1 \times
p}.\end{equation}
 \item[(iii)]
If  ${\rm rank}\,{\rm {\bf A}} = n$ and  ${\rm rank}\,{\rm {\bf
B}} = r_{2} < p$, then
\begin{equation}
\label{eq:AXB_detA*A_d^B} x_{ij}={\frac{{{ {{\rm{cdet}} _{i} \left( {\left( {{\rm {\bf A}}^{\sharp} {\rm {\bf A}}} \right)_{.\,i} \left( {{{\rm {\bf
d}}}\,_{.\,j}^{{\rm {\bf B}}}} \right)} \right)
} } }}{{{\rm det}( {\bf A}^{\sharp}  {\bf A})\cdot\sum\limits_{\alpha \in I_{r_{2},p}}{{\left| {\left( {{\rm
{\bf B}}{\rm {\bf B}}^{\sharp} } \right) _{\alpha} ^{\alpha} }
\right|}}} }},
\end{equation}
or
\begin{equation}\label{AXB_detA*A_d^A}
x_{ij}={\frac{{{\sum\limits_{\alpha \in I_{r_{2},p} {\left\{ {j}
\right\}}} {\rm rdet}_{j}{ \left( {\left( {{\rm {\bf B}}{\rm {\bf
B}}^{\sharp} } \right)_{\,j\,.} \left( {{{\rm {\bf d}}}\,_{i\,.}^{{\rm
{\bf A}}}} \right)} \right)\,_{\alpha} ^{\alpha} } }}}{{{\rm det}( {\bf A}^{\sharp}  {\bf A})\cdot\sum\limits_{\alpha \in I_{r_{2},p}}{{\left| {\left(
{{\rm {\bf B}}{\rm {\bf B}}^{\sharp} } \right) _{\alpha} ^{\alpha} }
\right|}}} }},
\end{equation}
where  $ {{\rm {\bf d}}_{.\,j}^{{\rm {\bf B}}}}$ is
(\ref{eq:def_d^B_m})
   and
  ${\rm {\bf
d}}_{i\,.}^{{\rm {\bf A}}}$ is (\ref{eq:def_d^A}).

\item[(iv)]
If ${\rm rank}\,{\rm {\bf A}} = r_{1} <  n$ and ${\rm rank}\,{\rm
{\bf B}} =  p$, then
\begin{equation}
\label{eq:AXB_detBB*_d^A} x_{i\,j} = {\frac{{{\rm rdet} _{j} ({\rm
{\bf B}}{\rm {\bf B}}^{\sharp} )_{j.\,} \left( {{\rm {\bf
d}}_{i\,.}^{{\rm {\bf A}}}} \right)}}{{{\sum\limits_{\beta
\in J_{r_{1},n}} {{\left| {\left( {{\rm {\bf A}}^{\sharp} {\rm {\bf
A}}} \right) _{\beta} ^{\beta} } \right|}}\cdot  {\rm det} (
{\bf B} {\bf B}^{\sharp} )}}}},
\end{equation}
or
\begin{equation} \label{eq:AXB_detBB*_d^B} x_{i\,j}=
{\frac{{{\sum\limits_{\beta \in J_{r_{1},\,n} {\left\{
{i} \right\}}} {{\rm{cdet}} _{i} \left( {\left( {{\rm {\bf A}}^{\sharp} {\rm {\bf A}}} \right)_{\,.\,i} \left( {{{\rm {\bf
d}}}\,_{.\,j}^{{\rm {\bf B}}}} \right)} \right) _{\beta} ^{\beta}
} } }}{{{\sum\limits_{\beta
\in J_{r_{1},n}} {{\left| {\left( {{\rm {\bf A}}^{\sharp} {\rm {\bf
A}}} \right)_{\beta} ^{\beta} } \right|}}\cdot {\rm det} (
{\bf B} {\bf B}^{\sharp} )}}}},
\end{equation}
 \noindent where  $ {{\rm {\bf d}}_{.\,j}^{{\rm {\bf B}}}}$ is
(\ref{eq:d_B_p})
   and
 ${\rm {\bf
d}}_{i\,.}^{{\rm {\bf A}}}$ is (\ref{eq:d_A<n}).

\end{enumerate}
\end{theorem}
{\bf Proof.}(i) If ${\rm {\bf A}} \in {\rm {\mathbb{H}}}_{r_{1}}^{m\times n}
$, ${\rm {\bf B}} \in {\rm {\mathbb{H}}}_{r_{2}}^{p\times q} $ and
$ r_{1} <  n$, $ r_{2} < p$, then, by Theorem
\ref{th:det_rep_A_mn}, the weighted  Moore-Penrose inverses ${\rm {\bf
A}}^{\dag} = \left( {{a}_{ij}^{ \ddag} } \right) \in
{\mathbb{H}}^{n\times m} $ and ${\rm {\bf B}}^{\dag} = \left(
{{b}_{ij}^{ \ddag} } \right) \in {\rm {\mathbb{H}}}^{q\times p} $ possess
the following determinantal representations, respectively,
\begin{gather}\label{eq:a+nf}
 {a}_{ij}^{\ddag}  = {\frac{{{\sum\limits_{\beta
\in J_{r_{1},\,n} {\left\{ {i} \right\}}} {{\rm{cdet}} _{i} \left(
{\left( {{\rm {\bf A}}^{\sharp} {\rm {\bf A}}} \right)_{\,. \,i}
\left( {{\rm {\bf a}}_{.j}^{\sharp} }  \right)} \right){\kern 1pt}
{\kern 1pt} _{\beta} ^{\beta} } } }}{{{\sum\limits_{\beta \in
J_{r_{1},\,n}} {{\left| {\left( {{\rm {\bf A}}^{\sharp} {\rm {\bf A}}}
\right){\kern 1pt} _{\beta} ^{\beta} }  \right|}}} }}},
\\
\label{eq:b+nf}
{b}_{ij}^{\ddag}  =
{\frac{{{\sum\limits_{\alpha \in I_{r_{2},p} {\left\{ {j}
\right\}}} {{\rm{rdet}} _{j} \left( {({\rm {\bf B}}{\rm {\bf B}}^{\sharp} )_{j\,.\,} ({\rm {\bf b}}_{i.\,}^{\sharp} )} \right)\,_{\alpha}
^{\alpha} } }}}{{{\sum\limits_{\alpha \in I_{r_{2},p}}  {{\left|
{\left( {{\rm {\bf B}}{\rm {\bf B}}^{\sharp} } \right){\kern 1pt}
_{\alpha} ^{\alpha} } \right|}}} }}}.
\end{gather}
By Lemma \ref{lem:sol_res},  ${\bf X}= {\bf A}^{\dag}_{M,N}{\bf D}{\bf B}^{\dag}_{P,Q}$ and entries of ${\rm {\bf
X}}=(x_{ij})$ are

\begin{equation}
\label{eq:sum+} x_{ij} = {{\sum\limits_{s = 1}^{q} {\left(
{{\sum\limits_{k = 1}^{m} { {a}_{ik}^{\ddag} d_{ks}} } } \right)}}
 {b}_{sj}^{\ddag}}.
\end{equation}
 for all $i=1,...,n$, $j=1,...,p$.

Denote  by ${\hat {\bf d}}_{.s}$ the $s$th column of ${\rm
{\bf A}}^{ \sharp}{\rm {\bf D}}=:\hat{{\rm {\bf D}}}=
(\hat{d}_{ij})\in {\mathbb{H}}^{n\times q}$ for all $s=1,...,q$. It follows from ${\sum\limits_{k} { {\rm {\bf a}}_{.\,k}^{\sharp}}d_{ks} }={\hat {\bf d}}_{.\,s}$ that

{\footnotesize{\begin{multline}\label{eq:sum_cdet}
\sum\limits_{k = 1}^{m} {{a}_{ik}^{\ddag} d_{ks}}=\sum\limits_{k =
1}^{m}{\frac{{{\sum\limits_{\beta \in J_{r_{1},\,n} {\left\{ {i}
\right\}}} {{\rm{cdet}} _{i} \left( {\left( {{\rm {\bf A}}^{\sharp}
{\rm {\bf A}}} \right)_{. \,i} \left( {{\rm {\bf a}}_{.k}^{\sharp} }
\right)} \right) _{\beta} ^{\beta} } }
}}{{{\sum\limits_{\beta \in J_{r_{1},\,n}} {{\left| {\left( {{\rm
{\bf A}}^{\sharp} {\rm {\bf A}}} \right) _{\beta} ^{\beta}
}  \right|}}} }}}\cdot d_{ks}=
\\
{\frac{{{\sum\limits_{\beta \in J_{r_{1},\,n} {\left\{ {i}
\right\}}}\sum\limits_{k = 1}^{m} {{\rm{cdet}} _{i} \left( {\left(
{{\rm {\bf A}}^{\sharp} {\rm {\bf A}}} \right)_{\,. \,i} \left( {{\rm
{\bf a}}_{.k}^{\sharp} } \right)} \right)
_{\beta} ^{\beta} } } }\cdot d_{ks}}{{{\sum\limits_{\beta \in
J_{r_{1},\,n}} {{\left| {\left( {{\rm {\bf A}}^{\sharp} {\rm {\bf A}}}
\right) _{\beta} ^{\beta} }  \right|}}}
}}}={\frac{{{\sum\limits_{\beta \in J_{r_{1},\,n} {\left\{ {i}
\right\}}} {{\rm{cdet}} _{i} \left( {\left( {{\rm {\bf A}}^{\sharp}
{\rm {\bf A}}} \right)_{\,. \,i} \left( {{\hat {\bf d}}_{.\,s}}
\right)} \right) _{\beta} ^{\beta} } }
}}{{{\sum\limits_{\beta \in J_{r_{1},\,n}} {{\left| {\left( {{\rm
{\bf A}}^{\sharp} {\rm {\bf A}}} \right) _{\beta} ^{\beta}
}  \right|}}} }}}.
\end{multline}}}
Suppose ${\rm {\bf e}}_{s.}$ and ${\rm {\bf e}}_{.\,s}$ are
 the unit row-vector and the unit column-vector, respectively, such that all their
components are $0$, except the $s$th components, which are $1$.
Substituting  (\ref{eq:sum_cdet}) and (\ref{eq:b+nf}) in
(\ref{eq:sum+}), we obtain

{\footnotesize{\[
x_{ij} =\sum\limits_{s = 1}^{q}{\frac{{{\sum\limits_{\beta \in
J_{r_{1},\,n} {\left\{ {i} \right\}}} {{\rm{cdet}} _{i} \left(
{\left( {{\rm {\bf A}}^{\sharp} {\rm {\bf A}}} \right)_{\,. \,i}
\left( {{\hat {\bf d}}_{.\,s}} \right)} \right) _{\beta} ^{\beta} } } }}{{{\sum\limits_{\beta \in
J_{r_{1},\,n}} {{\left| {\left( {{\rm {\bf A}}^{\sharp} {\rm {\bf A}}}
\right) _{\beta} ^{\beta} }  \right|}}}
}}}{\frac{{{\sum\limits_{\alpha \in I_{r_{2},p} {\left\{ {j}
\right\}}} {{\rm{rdet}} _{j} \left( {({\rm {\bf B}}{\rm {\bf B}}^{\sharp} )_{j\,.\,} ({\rm {\bf b}}_{s.\,}^{\sharp} )} \right)_{\alpha}
^{\alpha} } }}}{{{\sum\limits_{\alpha \in I_{r_{2},p}}  {{\left|
{\left( {{\rm {\bf B}}{\rm {\bf B}}^{\sharp} } \right)
_{\alpha} ^{\alpha} } \right|}}} }}}.
\]}}
Since \begin{equation*}{\hat{\bf
d}}_{.\,s}=\sum\limits_{l = 1}^{n}{\rm {\bf e}}_{.\,l}\hat{
d_{ls}},\,  {\rm {\bf b}}_{s.\,}^{\sharp}=\sum\limits_{t =
1}^{p}b_{st}^{\sharp}{\rm {\bf
e}}_{t.},\,\sum\limits_{s=1}^{q}{\hat d}_{ls}b_{st}^{\sharp}=\widetilde{d}_{lt},\end{equation*}
then we have

{\footnotesize{\[
x_{ij} = \]
\[{\frac{{ \sum\limits_{s = 1}^{q}\sum\limits_{t =
1}^{p} \sum\limits_{l = 1}^{n} {\sum\limits_{\beta \in
J_{r_{1},\,n} {\left\{ {i} \right\}}} {{\rm{cdet}} _{i} \left(
{\left( {{\rm {\bf A}}^{\sharp} {\rm {\bf A}}} \right)_{\,. \,i}
\left( {\rm {\bf e}}_{.\,l} \right)} \right)  _{\beta} ^{\beta} } } }{\kern 1pt}{\hat
d}_{ls}b_{st}^{\sharp}{\sum\limits_{\alpha \in I_{r_{2},p} {\left\{ {j}
\right\}}} {{\rm{rdet}} _{j} \left( {({\rm {\bf B}}{\rm {\bf B}}^{\sharp} )_{j\,.\,} ({\rm {\bf e}}_{t.} )} \right)_{\alpha} ^{\alpha}
} } }{{{\sum\limits_{\beta \in J_{r_{1},\,n}} {{\left| {\left(
{{\rm {\bf A}}^{\sharp} {\rm {\bf A}}} \right) _{\beta}
^{\beta} }  \right|}}} }{{\sum\limits_{\alpha \in I_{r_{2},p}}
{{\left| {\left( {{\rm {\bf B}}{\rm {\bf B}}^{\sharp} } \right) _{\alpha} ^{\alpha} } \right|}}} }}    }=
\]}}
\begin{equation}\label{eq:x_ij}
{\frac{{ \sum\limits_{t = 1}^{p} \sum\limits_{l = 1}^{n}
{\sum\limits_{\beta \in J_{r_{1},\,n} {\left\{ {i} \right\}}}
{{\rm{cdet}} _{i} \left( {\left( {{\rm {\bf A}}^{\sharp} {\rm {\bf
A}}} \right)_{\,. \,i} \left( {\rm {\bf e}}_{.\,l} \right)}
\right) _{\beta} ^{\beta} } }
}\,\,\widetilde{d}_{lt}{\sum\limits_{\alpha \in I_{r_{2},p}
{\left\{ {j} \right\}}} {{\rm{rdet}} _{j} \left( {({\rm {\bf
B}}{\rm {\bf B}}^{\sharp} )_{j\,.\,} ({\rm {\bf e}}_{t.} )}
\right)_{\alpha} ^{\alpha} } } }{{{\sum\limits_{\beta \in
J_{r_{1},\,n}} {{\left| {\left( {{\rm {\bf A}}^{\sharp} {\rm {\bf A}}}
\right){\kern 1pt} _{\beta} ^{\beta} }  \right|}}}
}{{\sum\limits_{\alpha \in I_{r_{2},p}} {{\left| {\left( {{\rm
{\bf B}}{\rm {\bf B}}^{\sharp} } \right) _{\alpha}
^{\alpha} } \right|}}} }}    }.
\end{equation}
Denote by

{\footnotesize{\[
 d^{{\rm {\bf A}}}_{it}:= \]
\[
{\sum\limits_{\beta \in J_{r_{1},\,n} {\left\{ {i} \right\}}}
{{\rm{cdet}} _{i} \left( {\left( {{\rm {\bf A}}^{\sharp} {\rm {\bf
A}}} \right)_{. \,i} \left( \widetilde{{\rm {\bf d}}}_{.t}
\right)} \right) _{\beta} ^{\beta} } }= \sum\limits_{l
= 1}^{n} {\sum\limits_{\beta \in J_{r_{1},n} {\left\{ {i}
\right\}}} {{\rm{cdet}} _{i} \left( {\left( {{\rm {\bf A}}^{\sharp}
{\rm {\bf A}}} \right)_{\,. \,i} \left( {\rm {\bf e}}_{.l}
\right)} \right) _{\beta} ^{\beta} } }
\widetilde{d}_{lt}
\]}}
the $t$th component  of a row-vector ${\rm {\bf d}}^{{\rm {\bf
A}}}_{i\,.}= (d^{{\rm {\bf A}}}_{i1},...,d^{{\rm {\bf A}}}_{ip})$
for all $t=1,...,p$. Substituting it in (\ref{eq:x_ij}),
we have
\[x_{ij} ={\frac{{ \sum\limits_{t = 1}^{p}
 d^{{\rm {\bf A}}}_{it}
}{\sum\limits_{\alpha \in I_{r_{2},p} {\left\{ {j} \right\}}}
{{\rm{rdet}} _{j} \left( {({\rm {\bf B}}{\rm {\bf B}}^{\sharp}
)_{j\,.\,} ({\rm {\bf e}}_{t.} )} \right)\,_{\alpha} ^{\alpha} } }
}{{{\sum\limits_{\beta \in J_{r_{1},\,n}} {{\left| {\left( {{\rm
{\bf A}}^{\sharp} {\rm {\bf A}}} \right){\kern 1pt} _{\beta} ^{\beta}
}  \right|}}} }{{\sum\limits_{\alpha \in I_{r_{2},p}} {{\left|
{\left( {{\rm {\bf B}}{\rm {\bf B}}^{\sharp} } \right){\kern 1pt}
_{\alpha} ^{\alpha} } \right|}}} }}    }.
\]
Since $\sum\limits_{t = 1}^{p}
 d^{{\rm {\bf A}}}_{it}{\rm {\bf e}}_{t.}={\rm {\bf
d}}^{{\rm {\bf A}}}_{i\,.}$, then it follows (\ref{eq:d^A}).

If we denote by

{\footnotesize{\begin{equation*}
 d^{{\rm {\bf B}}}_{lj}:=
\sum\limits_{t = 1}^{p}\widetilde{d}_{lt}{\sum\limits_{\alpha \in
I_{r_{2},p} {\left\{ {j} \right\}}} {{\rm{rdet}} _{j} \left(
{({\rm {\bf B}}{\rm {\bf B}}^{\sharp} )_{j\,.\,} ({\rm {\bf e}}_{t.}
)} \right)_{\alpha} ^{\alpha} } }={\sum\limits_{\alpha \in
I_{r_{2},p} {\left\{ {j} \right\}}} {{\rm{rdet}} _{j} \left(
{({\rm {\bf B}}{\rm {\bf B}}^{\sharp} )_{j\,.\,} (\widetilde{{\rm {\bf
d}}}_{l.} )} \right)_{\alpha} ^{\alpha} } }
\end{equation*}}}
\noindent the $l$th component  of a column-vector ${\rm {\bf
d}}^{{\rm {\bf B}}}_{.\,j}= (d^{{\rm {\bf B}}}_{1j},...,d^{{\rm
{\bf B}}}_{jn})^{T}$ for all $l=1,...,n$ and substitute it
in (\ref{eq:x_ij}), we obtain
\[x_{ij} ={\frac{{  \sum\limits_{l = 1}^{n}
{\sum\limits_{\beta \in J_{r_{1},\,n} {\left\{ {i} \right\}}}
{{\rm{cdet}} _{i} \left( {\left( {{\rm {\bf A}}^{\sharp} {\rm {\bf
A}}} \right)_{\,. \,i} \left( {\rm {\bf e}}_{.\,l} \right)}
\right){\kern 1pt}  _{\beta} ^{\beta} } } }\,\,d^{{\rm {\bf
B}}}_{lj} }{{{\sum\limits_{\beta \in J_{r_{1},\,n}} {{\left|
{\left( {{\rm {\bf A}}^{\sharp} {\rm {\bf A}}} \right){\kern 1pt}
_{\beta} ^{\beta} }  \right|}}} }{{\sum\limits_{\alpha \in
I_{r_{2},p}} {{\left| {\left( {{\rm {\bf B}}{\rm {\bf B}}^{\sharp} }
\right){\kern 1pt} _{\alpha} ^{\alpha} } \right|}}} }}    }.
\]
Since $\sum\limits_{l = 1}^{n}{\rm {\bf e}}_{.l}
 d^{{\rm {\bf B}}}_{lj}={\rm {\bf
d}}^{{\rm {\bf B}}}_{.\,j}$, then it follows (\ref{eq:d^B}).

(ii) If ${\rm rank}\,{\rm {\bf A}} = n$ and ${\rm rank}\,{\rm {\bf
B}} = p$,
  then by Theorem \ref{th:det_rep_A_mn} the weighted  Moore-Penrose inverses ${\rm {\bf
A}}_{M,N}^{ \dag} = \left( {{a}_{ij}^{ \ddag} } \right) \in
{\mathbb{H}}^{n\times m} $ and ${\rm {\bf B}}_{P,Q}^{\dag} = \left(
{{b}_{ij}^{\ddag} } \right) \in {\rm {\mathbb{H}}}^{q\times p} $ possess
the following determinantal representations, respectively,
\begin{gather}\label{eq:A1A_f}
 {{a}_{ij}^{\ddag} }   = {\frac{{{\rm{cdet}} _{i} ({{\rm {\bf A}}^{\sharp} {\rm {\bf
A}}})_{.\,i} \left( {{\rm {\bf a}}_{.\,j}^{\sharp} }  \right)} }{{ {\rm det} ({{\rm {\bf A}}^{\sharp} {\rm {\bf
A}}})}}}\\\label{eq:BB1_f}
  {{b}_{ij}^{\ddag} }   =    {\frac{{{\rm{rdet}} _{j} ({\rm {\bf B}} {\rm {\bf B}}^{\sharp})_{j.} \left( {{\rm {\bf b}}_{i.}^{\sharp} }  \right)}}{{{\rm det} {( {\bf B} {\bf B}^{\sharp}) }}}}.
\end{gather}
By their substituting in (\ref{eq:sum+}) and
pondering ahead as in the previous case, we obtain (\ref{eq:AXB_cdetA*A}) and (\ref{eq:AXB_rdetBB*}).

(iii)
 If ${\rm {\bf A}} \in {\rm {\mathbb{H}}}_{r_{1}}^{m\times n}
$, ${\rm {\bf B}} \in {\rm {\mathbb{H}}}_{r_{2}}^{p\times q} $ and
$ r_{1} =n$, $ r_{2} < p$, then, for the weighted  Moore-Penrose inverses ${\bf A}_{M,N}^{\dag}$ and ${\bf B}_{P,Q}^{\dag}$, the determinantal representations (\ref{eq:A1A_f})
 and  (\ref{eq:a+nf})  are more applicable to use, respectively.
By their substituting in (\ref{eq:sum+}) and
pondering ahead as in the previous case, we finally obtain (\ref{eq:AXB_detA*A_d^B}) and (\ref{AXB_detA*A_d^A}) as well.

(iv) In this case  for  ${\bf A}_{M,N}^{\dag}$ and ${\bf B}_{P,Q}^{\dag}$, we use the determinantal representations
(\ref{eq:A1A_f}) and (\ref{eq:b+nf}), respectively.
$\Box$

\begin{corollary}\label{cor:AX} Suppose that ${\bf A}\in {\mathbb{H}}^{m\times n}_{{r}_1}$, ${\bf D}\in {\mathbb{H}}^{m\times p}$,   ${\bf M}$, ${\bf N}$ are Hermitian positive definite matrices of order $m$  and $n$, respectively, ${\bf A}^{\sharp}{\bf A}$ is Hermitian. Denote
${ \widehat{\bf{D}}}=  {\bf
A}^{\sharp} {\bf D}$. If ${\bf D}\subset  \mathcal{R}_{r}({\bf A}{\bf A}^{\sharp})$ and ${\bf D}\subset  \mathcal{R}_{l}({\bf A}^{\sharp}{\bf A})$,
\begin{gather}\label{eq_AX}
{\bf A}{\bf X}={\bf D},\\\label{res_AX}
\mathcal{R}_{r}({\bf X})\subset {\bf N}^{-1}\mathcal{R}_{r}({\bf A}^{*}),\,\mathcal{R}_{l}({\bf X})\subset \mathcal{R}_{l}({\bf A}^{*}){\bf M},
\end{gather}
then the unique solution of (\ref{eq_AX}) with the restrictions (\ref{res_AX}) is
\begin{equation*}
{\bf X}= {\bf A}^{\dag}_{M,N}{\bf D}\end{equation*}
which possess the following determinantal representations.
\begin{enumerate}
\item[(i)] If ${\rm rank}\,{\rm {\bf A}} = r_{1} <  n$, then
\begin{equation*}
x_{ij} = {\frac{{{\sum\limits_{\beta \in J_{{{r}_1},\,n} {\left\{
{i} \right\}}} {{\rm{cdet}} _{i} \left( {\left( {{\rm {\bf A}}^{\sharp} {\rm {\bf A}}} \right)_{.\,i} \left( {{ \widehat{{\bf
d}}}\,_{.\,j}} \right)} \right) _{\beta} ^{\beta}
} } }}{{{\sum\limits_{\beta \in J_{{{r}_1},n}} {{\left| {\left(
{{\rm {\bf A}}^{\sharp} {\rm {\bf A}}} \right)_{\beta} ^{\beta} }
\right|}}}} }},
\end{equation*}
where
${\widehat {\bf d}}_{.j}$ are  the $j$th
column  of ${ \widehat {\bf {D}}}$ for all $i =1,...,n $, $j =1,...,p $.
\item[(ii)] If ${\rm rank}\,{\rm {\bf A}} = n$, then
\begin{equation*}
 x_{i\,j} = \frac{{{\rm cdet} _{i} (
{\bf A}^{\sharp}  {\bf A})_{.\,i\,} \left( {{\widehat {\bf
d}}_{.j}} \right)}}{{\rm det} ({\bf A}^{\sharp}  {\bf A})},
\end{equation*}
\end{enumerate}
\end{corollary}
{\bf Proof.} The proof follows evidently from Theorem \ref{th:cr_rul_AXB}  when ${\bf B}$ be removed,  and unit matrices  insert instead ${\bf P}$, ${\bf Q}$.

\begin{corollary}\label{cor:XB} Suppose that ${\bf B}\in {\mathbb{H}}^{p\times q}_{{r}_2}$, ${\bf D}\in {\mathbb{H}}^{n\times q}$, ${\bf P}$, and ${\bf Q}$ are Hermitian positive definite matrices of order  $p$ and $q$, respectively, ${\bf B}{\bf B}^{\sharp}$ is Hermitian. Denote
${ \check{\bf{D}}}=  {\bf D}{\bf B}^{\sharp}$. If ${\bf D}\subset  \mathcal{R}_{r}({\bf B}^{\sharp}{\bf B})$ and ${\bf D}\subset  \mathcal{R}_{l}({\bf B}{\bf B}^{\sharp})$,
\begin{gather}\label{eq_XB}
{\bf X}{\bf B}={\bf D},\\\label{res_XB}
\mathcal{N}_{r}({\bf X})\supset {\bf P}^{-1}\mathcal{N}_{r}({\bf B}^{*}),\,\mathcal{N}_{l}({\bf X})\supset \mathcal{N}_{l}({\bf B}^{*}){\bf Q} ,
\end{gather}
then the unique solution of (\ref{eq_XB}) with the restrictions (\ref{res_XB}) is
\begin{equation*}
{\bf X}= {\bf D}{\bf B}^{\dag}_{P,Q}\end{equation*}
which possess the following determinantal representations.
\begin{enumerate}
\item[(i)] If ${\rm rank}\,{\bf B} = r_{2} <  p$, then
\begin{equation*}
 x_{ij}={\frac{{{\sum\limits_{\alpha
\in I_{r_{2},q} {\left\{ {j} \right\}}} {{\rm{rdet}} _{j} \left(
{\left( {{\rm {\bf B}}{\rm {\bf B}}^{\sharp} } \right)_{j\,.} \left(
{{{\check {\bf d}}}\,_{i\,.}} \right)}
\right)\,_{\alpha} ^{\alpha} } }}}{{{\sum\limits_{\alpha \in
I_{r_{2},q}} {{\left| {\left( {{\rm {\bf B}}{\rm {\bf B}}^{\sharp} }
\right) _{\alpha} ^{\alpha} } \right|}}} }}},
\end{equation*}
where
${\check {\bf d}}_{i.}$ are  the $i$th
row  of ${ \check {\bf {D}}}$ for all $i =1,...,n $, $j =1,...,p $.
\item[(ii)] If ${\rm rank}\,{\bf B} =   p$, then
\begin{equation*}
 x_{i\,j} = \frac{{{\rm rdet} _{j} \left( {{\rm {\bf B}}{\rm {\bf B}}^{\sharp} }\right)_{j.} \left( {{\check {\bf
d}}_{i.}} \right)}}{{\rm det} \left( {{\rm {\bf B}}{\rm {\bf B}}^{\sharp} }
\right)}.
\end{equation*}
\end{enumerate}
\end{corollary}
{\bf Proof.} The proof follows evidently from Theorem \ref{th:cr_rul_AXB} when  ${\bf A}$  be removed and unit matrices  insert instead ${\bf M}$, ${\bf N}$.

\subsection{The Case of Both Non-Hermitian Matrices ${\bf A}^{\sharp}{\bf A}$ and ${\bf B}{\bf B}^{\sharp}$.}
Denote  $\widetilde{{\bf A}}:=
{\bf M}^{\frac{1}{2}}{\bf A}{\bf N}^{-\frac{1}{2}}
 =\left(\widetilde{a}_{ij}\right)\in {\mathbb{H}}^{m\times n} $,  $\widetilde{{\bf A}}^{*}={\bf N}^{-\frac{1}{2}}{\bf A}^{*}{\bf M}^{\frac{1}{2}}$ and  $\widetilde{{\bf B}}:=
{\bf P}^{\frac{1}{2}}{\bf B}{\bf Q}^{-\frac{1}{2}}
 =\left(\widetilde{a}_{ij}\right)\in {\mathbb{H}}^{p\times q} $,  $\widetilde{{\bf B}}^{*}={\bf Q}^{-\frac{1}{2}}{\bf B}^{*}{\bf P}^{\frac{1}{2}}$.
\begin{theorem}\label{th:cr_rul_AXB_n} Let ${\bf A}^{\sharp}{\bf A}$ and ${\bf B}{\bf B}^{\sharp}$ be both non-Hermitian.  Then the  solution  (\ref{eq:sol_res})  possess the following determinantal representations.
\begin{enumerate}
\item[(i)] If ${\rm rank}\,{\rm {\bf A}} = r_{1} <  n$ and
 ${\rm rank}\,{\rm {\bf B}} = r_{2} < p$, then
\begin{equation}\label{eq:d^B_n}
x_{ij} = {\frac{{    {\sum\limits_{k}n^{(-\frac{1}{2})}_{ik}\sum\limits_{\beta \in
J_{r_{1},\,n} {\left\{ {k} \right\}}} {\rm{cdet}} _{k} \left( {\left(
\widetilde{{\bf A}}^{*}\widetilde{{\bf A}}
 \right)_{\,. \,k} \left({{\rm {\bf d}}}^{{\rm {\bf
B}}}_{.j}
\right)} \right) _{\beta} ^{\beta}} } }{{{\sum\limits_{\beta \in J_{r_{1},\,n}} {{\left|
{\left(
\widetilde{{\bf A}}^{*}\widetilde{{\bf A}}
 \right){\kern 1pt}
_{\beta} ^{\beta} }  \right|}}} }{{\sum\limits_{\alpha \in
I_{r_{2},p}} {{\left| {\left( \widetilde{{\bf B}}\widetilde{{\bf B}}^{*} \right) _{\alpha} ^{\alpha} } \right|}}} }}    },
\end{equation} or

\begin{equation}\label{eq:d^A_n}
 x_{ij}={\frac{{\sum\limits_{l}{{\sum\limits_{\alpha \in
I_{r_{2},\,p} {\left\{ {l} \right\}}} {{\rm{rdet}} _{l} \left( {\left(
\widetilde{{\bf B}}\widetilde{{\bf B}}^{*}
 \right)_{l. \,}  ({{\bf d}}^{{\rm {\bf A}}}_{i.})  } \right)
_{\alpha} ^{\alpha} } } }\cdot m^{(\frac{1}{2})}_{lj}}
}{{{\sum\limits_{\beta \in
J_{r_{1},\,n}} {{\left| {\left(
\widetilde{{\bf A}}^{*}\widetilde{{\bf A}}
 \right)  _{\beta} ^{\beta} }  \right|}}}
}{{\sum\limits_{\alpha \in I_{r_{2},p}} {{\left| {\left( \widetilde{{\bf B}}\widetilde{{\bf B}}^{*} \right) _{\alpha} ^{\alpha} }  \right|}}} }}     },
\end{equation}

where
\begin{equation} \label{eq:def_d^B_m_n}
   {{{\rm {\bf d}}}_{.\,j}^{{\rm {\bf B}}}}= \left( {\sum\limits_{l}{{\sum\limits_{\alpha \in
I_{r_{2},\,p} {\left\{ {l} \right\}}} {{\rm{rdet}} _{l} \left( {\left(
\widetilde{{\bf B}}\widetilde{{\bf B}}^{*}
 \right)_{l. \,}  ({\widetilde {\bf d}}_{t.})  } \right)
_{\alpha} ^{\alpha} } } }\cdot m^{(\frac{1}{2})}_{lj}}
\right)\in{\rm {\mathbb{H}}}^{n \times
1}\end{equation}
\begin{equation} \label{eq:d_A<n_n}  {{{\rm {\bf d}}}_{i\,.}^{{\rm {\bf A}}}}=\left(
 {\sum\limits_{k}n^{(-\frac{1}{2})}_{ik}\sum\limits_{\beta \in
J_{r_{1},\,n} {\left\{ {k} \right\}}} {\rm{cdet}} _{k} \left( {\left(
\widetilde{{\bf A}}^{*}\widetilde{{\bf A}}
 \right)_{\,. \,k} \left({{\widetilde {\bf d}}}_{.f}
\right)} \right) _{\beta} ^{\beta}}
 \right)\in{\rm {\mathbb{H}}}^{1 \times
p}
\end{equation}
are the column-vector and the row-vector, respectively.  ${\tilde{{\rm {\bf
 d}}}_{t.}}$ and
${\tilde{{\rm {\bf d}}}_{.f}}$ are the $t$th row   and the $f$th
column  of $\widetilde{{\bf D}}:={\bf N}^{-\frac{1}{2}}{\bf A}^{*}{\bf M}{\bf D}{\bf Q}^{-1}{\bf B}^{*}{\bf P}^{\frac{1}{2}}=(\widetilde{d}_{ij})\in {\mathbb{H}}^{n\times p}$ for all $t =1,...,n $, $f =1,...,p $.

\item[(ii)] If ${\rm rank}\,{\rm {\bf A}} = n$ and ${\rm rank}\,{\rm {\bf B}} = p$, then
\begin{equation}
\label{eq:AXB_cdetA*A_n} x_{i\,j} = \frac{{{\rm cdet} _{i} ({\bf A}^{*}{\bf M}{\bf A})_{.\,i\,} \left( {{\rm {\bf
d}}_{.j}^{{\rm {\bf B}}}} \right)}}{{\rm det} ({\bf A}^{*}{\bf M}{\bf A})\cdot {\rm det} ({\bf B}{\bf Q}^{-1}{\bf B}^{*})},
\end{equation}
or
\begin{equation}
\label{eq:AXB_rdetBB*_n} x_{i\,j} = \frac{{\rm rdet} _{j} ({\bf B}{\bf Q}^{-1}{\bf B}^{*})_{j.\,} \left(  {\bf
d}_{i\,.}^{ {\bf A}} \right)}{{\rm det} ({\bf A}^{*}{\bf M}{\bf A})\cdot {\rm det} ({\bf B}{\bf Q}^{-1}{\bf B}^{*})},
\end{equation}
 \noindent where   \begin{equation}\label{eq:d_B_p_n} {\bf d}_{.j}^ {\bf B} : =
\left( {\rm
rdet} _{j} ({\bf B}{\bf Q}^{-1}{\bf B}^{*})_{j.\,} \left(
{\widetilde {\bf d}}_{t\,.} \right) \right)\in{\mathbb{H}}^{n \times
1},\end{equation}
\begin{equation}\label{eq:def_d^A_n}{\bf d}_{i\,.}^ {\bf A} : = \left(
{\rm cdet}
_{i} ({\bf A}^{*}{\bf M}{\bf A})_{.\,i} \left(
{\widetilde {\bf d}}_{.f} \right) \right)\in{\rm {\mathbb{H}}}^{1 \times
p},\end{equation}
 ${\widetilde{\bf d}}_{t\,.}$ and ${\widetilde{\bf d}}_{.f}$ are the $t$th row and the $f$th column of ${\widetilde {\bf D}}:={\bf A}^{*}{\bf M}{\bf D}{\bf Q}^{-1}{\bf B}^{*}\in {\mathbb{H}}^{n\times p}$, respectively.
 \item[(iii)]
If  ${\rm rank}\,{\rm {\bf A}} = n$ and  ${\rm rank}\,{\rm {\bf
B}} = r_{2} < p$, then
\begin{equation}
\label{eq:AXB_detA*A_d^B_n} x_{ij}={\frac{{{ {{\rm{cdet}} _{i} \left( {({\bf A}^{*}{\bf M}{\bf A})_{.\,i} \left( {{{\rm {\bf
d}}}\,_{.\,j}^{{\rm {\bf B}}}} \right)} \right)
} } }}{{{\rm det}({\bf A}^{*}{\bf M}{\bf A})\cdot\sum\limits_{\alpha \in I_{r_{2},p}}{{\left| {\left( \widetilde{{\bf B}}\widetilde{{\bf B}}^{*} \right) _{\alpha} ^{\alpha} }
\right|}}} }},
\end{equation}
or
\begin{equation}\label{AXB_detA*A_d^A_n}
x_{ij}={\frac{{\sum\limits_{l}{{\sum\limits_{\alpha \in
I_{r_{2},\,p} {\left\{ {l} \right\}}} {{\rm{rdet}} _{l} \left( {\left(
\widetilde{{\bf B}}\widetilde{{\bf B}}^{*}
 \right)_{l. \,}  ({{\bf d}}^{{\rm {\bf A}}}_{i.})  } \right)
_{\alpha} ^{\alpha} } } }\cdot m^{(\frac{1}{2})}_{lj}}
}{{{\rm det}({\bf A}^{*}{\bf M}{\bf A})\cdot\sum\limits_{\alpha \in I_{r_{2},p}}{{\left| {\left( \widetilde{{\bf B}}\widetilde{{\bf B}}^{*} \right) _{\alpha} ^{\alpha} }
\right|}}} }},
\end{equation}
where  $ {{\rm {\bf d}}_{.\,j}^{{\rm {\bf B}}}}$ is
(\ref{eq:def_d^B_m_n})
   and
  ${\rm {\bf
d}}_{i\,.}^{{\rm {\bf A}}}$ is (\ref{eq:def_d^A_n}).

\item[(iv)]
If ${\rm rank}\,{\rm {\bf A}} = r_{1} <  n$ and ${\rm rank}\,{\rm
{\bf B}} =  p$, then
\begin{equation}
\label{eq:AXB_detBB*_d^A_n} x_{i\,j} = {\frac{{{\rm rdet} _{j} ({\bf B}{\bf Q}^{-1}{\bf B}^{*})_{j.\,} \left( {{\rm {\bf
d}}_{i\,.}^{{\rm {\bf A}}}} \right)}}{{{\sum\limits_{\beta
\in J_{r_{1},n}} {{\left| {\left(
\widetilde{{\bf A}}^{*}\widetilde{{\bf A}}
 \right) _{\beta} ^{\beta} } \right|}}\cdot  {\rm det} ({\bf B}{\bf Q}^{-1}{\bf B}^{*})}}}},
\end{equation}
or
\begin{equation} \label{eq:AXB_detBB*_d^B_n} x_{i\,j}=
{\frac{{    {\sum\limits_{k}n^{(-\frac{1}{2})}_{ik}\sum\limits_{\beta \in
J_{r_{1},\,n} {\left\{ {k} \right\}}} {\rm{cdet}} _{k} \left( {\left(
\widetilde{{\bf A}}^{*}\widetilde{{\bf A}}
 \right)_{\,. \,k} \left({{\rm {\bf d}}}^{{\rm {\bf
B}}}_{.j}
\right)} \right) _{\beta} ^{\beta}} }}{{{\sum\limits_{\beta
\in J_{r_{1},n}} {{\left| {\left(
\widetilde{{\bf A}}^{*}\widetilde{{\bf A}}
 \right)_{\beta} ^{\beta} } \right|}}\cdot {\rm det} ({\bf B}{\bf Q}^{-1}{\bf B}^{*})}}}},
\end{equation}
 \noindent where  $ {{\rm {\bf d}}_{.\,j}^{{\rm {\bf B}}}}$ is
(\ref{eq:d_B_p_n})
   and
 ${\rm {\bf
d}}_{i\,.}^{{\rm {\bf A}}}$ is (\ref{eq:d_A<n_n}).

\end{enumerate}
\end{theorem}
{\bf Proof.}
(i) If ${\rm {\bf A}} \in {\rm {\mathbb{H}}}_{r_{1}}^{m\times n}
$, ${\rm {\bf B}} \in {\rm {\mathbb{H}}}_{r_{2}}^{p\times q} $ are both non-Hermitian, and
$ r_{1} <  n$, $ r_{2} < p$, then, by Theorem
\ref{th:det_rep_A_mn_non}, the weighted  Moore-Penrose inverses ${\rm {\bf
A}}^{\dag} = \left( {{a}_{ij}^{ \ddag} } \right) \in
{\mathbb{H}}^{n\times m} $ and ${\rm {\bf B}}^{\dag} = \left(
{{b}_{ij}^{ \ddag} } \right) \in {\rm {\mathbb{H}}}^{q\times p} $ posses
the following determinantal representations, respectively,
\begin{equation}\label{eq:a+nf_n}
 {{a}_{ij}^{\ddag} }=
{\frac{\sum\limits_{k}n^{(-\frac{1}{2})}_{ik}{{\sum\limits_{\beta \in
J_{r_{1},\,n} {\left\{ {k} \right\}}} {{\rm{cdet}} _{k} \left( {\left(
\widetilde{{\bf A}}^{*}\widetilde{{\bf A}}
 \right)_{\,. \,k}  \left(\widehat{{\bf a}}_{.j}\right)  } \right)
_{\beta} ^{\beta} } } }     }{{{\sum\limits_{\beta \in J_{r_{1},\,\,n}}
{{\left| {\left(
\widetilde{{\bf A}}^{*}\widetilde{{\bf A}}
 \right)  _{\beta} ^{\beta} }  \right|}}} }}},\end{equation}
 where $\widehat{{\bf a}}_{.j}$ is the $j$th column of ${\bf N}^{-\frac{1}{2}}{\bf A}^{*}{\bf M}$;
\begin{equation}\label{eq:b+nf_n}
 {{b}_{ij}^{\ddag} }=
 {\frac{\sum\limits_{l}{{\sum\limits_{\alpha \in
I_{r_{2},\,p} {\left\{ {l} \right\}}} {{\rm{rdet}} _{l} \left( {\left(
\widetilde{{\bf B}}\widetilde{{\bf B}}^{*}
 \right)_{l. \,}  (\widehat{{\bf b}}_{i.})  } \right)
_{\alpha} ^{\alpha} } } }\cdot m^{(\frac{1}{2})}_{lj}}{{{\sum\limits_{\alpha \in I_{_{2},\,\,p}}
{{\left| {\left( \widetilde{{\bf B}}\widetilde{{\bf B}}^{*} \right) _{\alpha} ^{\alpha} }  \right|}}} }}},
\end{equation}
where $\widehat{{\bf b}}_{i.}$ is the $i$th row of ${\bf Q}^{-1}{\bf B}^{*}{\bf P}^{\frac{1}{2}}$.
By Lemma \ref{lem:sol_res},  ${\bf X}= {\bf A}^{\dag}_{M,N}{\bf D}{\bf B}^{\dag}_{P,Q}$ and entries of ${\rm {\bf
X}}=(x_{ij})$ are

\begin{equation}
\label{eq:sum+_n} x_{ij} = {{\sum\limits_{s = 1}^{q} {\left(
{{\sum\limits_{t = 1}^{m} { {a}_{it}^{\ddag} d_{ts}} } } \right)}}
 {b}_{sj}^{\ddag}}.
\end{equation}
 for all $i=1,...,n$, $j=1,...,p$.

Denote  by ${\widehat{\bf d}}_{.s}$ the $s$th column of ${\bf N}^{-\frac{1}{2}}{\bf A}^{*}{\bf M}{\rm {\bf D}}=:\widehat{{\rm {\bf D}}}=
(\widehat{d}_{ij})\in {\mathbb{H}}^{n\times q}$ for all $s=1,...,q$. It follows from ${\sum\limits_{t} { {\widehat {\bf a}}_{.\,t}}d_{ts} }={\widehat {\bf d}}_{.\,s}$ that
\begin{multline}\label{eq:sum_cdet_n}
\sum\limits_{t = 1}^{m} { {a}_{it}^{\ddag} d_{ts}}=\sum\limits_{t =
1}^{m}{\frac{\sum\limits_{k}n^{(-\frac{1}{2})}_{ik}{{\sum\limits_{\beta \in
J_{r_{1},\,n} {\left\{ {i} \right\}}} {{\rm{cdet}} _{k} \left( {\left(
\widetilde{{\bf A}}^{*}\widetilde{{\bf A}}
 \right)_{\,. \,k}  \left(\widehat{{\bf a}}_{.t}\right)  } \right)
_{\beta} ^{\beta} } } }     }{{{\sum\limits_{\beta \in J_{r_{1},\,\,n}}
{{\left| {\left(
\widetilde{{\bf A}}^{*}\widetilde{{\bf A}}
 \right)  _{\beta} ^{\beta} }  \right|}}} }}}\cdot d_{ts}=
\\
{\frac{\sum\limits_{k}n^{(-\frac{1}{2})}_{ik}{{\sum\limits_{\beta \in
J_{r_{1},\,n} {\left\{ {i} \right\}}} {{\rm{cdet}} _{k} \left( {\left(
\widetilde{{\bf A}}^{*}\widetilde{{\bf A}}
 \right)_{\,. \,k}  \left(\widehat{{\bf d}}_{.s}\right)  } \right)
_{\beta} ^{\beta} } } }     }{{{\sum\limits_{\beta \in J_{r_{1},\,\,n}}
{{\left| {\left(
\widetilde{{\bf A}}^{*}\widetilde{{\bf A}}
 \right)  _{\beta} ^{\beta} }  \right|}}} }}}.
\end{multline}
Suppose ${\rm {\bf e}}_{s.}$ and ${\rm {\bf e}}_{.\,s}$ are
 the unit row-vector and the unit column-vector, respectively, such that all their
components are $0$, except the $s$th components, which are $1$.
Substituting  (\ref{eq:sum_cdet_n}) and (\ref{eq:b+nf_n}) in
(\ref{eq:sum+_n}), we obtain
\begin{multline*}
x_{ij} =\sum\limits_{s = 1}^{q}{\frac{\sum\limits_{k}n^{(-\frac{1}{2})}_{ik}{{\sum\limits_{\beta \in
J_{r_{1},\,n} {\left\{ {k} \right\}}} {{\rm{cdet}} _{k} \left( {\left(
\widetilde{{\bf A}}^{*}\widetilde{{\bf A}}
 \right)_{\,. \,k}  \left(\widehat{{\bf d}}_{.s}\right)  } \right)
_{\beta} ^{\beta} } } }     }{{{\sum\limits_{\beta \in J_{r_{1},\,\,n}}
{{\left| {\left(
\widetilde{{\bf A}}^{*}\widetilde{{\bf A}}
 \right)  _{\beta} ^{\beta} }  \right|}}} }}}\times\\{\frac{\sum\limits_{l}{{\sum\limits_{\alpha \in
I_{r_{2},\,p} {\left\{ {l} \right\}}} {{\rm{rdet}} _{l} \left( {\left(
\widetilde{{\bf B}}\widetilde{{\bf B}}^{*}
 \right)_{l. \,}  (\widehat{{\bf b}}_{s.})  } \right)
_{\alpha} ^{\alpha} } } }\cdot m^{(\frac{1}{2})}_{lj}}{{{\sum\limits_{\alpha \in I_{_{2},\,\,p}}
{{\left| {\left( \widetilde{{\bf B}}\widetilde{{\bf B}}^{*} \right) _{\alpha} ^{\alpha} }  \right|}}} }}}.
\end{multline*}
Since \begin{equation*}{\widehat{\bf
d}}_{.\,s}=\sum\limits_{l = 1}^{n}{\rm {\bf e}}_{.\,l}\widehat{
d}_{ls},\,  {\widehat {\bf b}}_{s.\,}=\sum\limits_{t =
1}^{p}\widehat{b}_{st}{\rm {\bf
e}}_{t.},\,\sum\limits_{s=1}^{q}{\widehat d}_{ls}\widehat{b}_{st}=\widetilde{d}_{lt},\end{equation*}
then we have
\begin{multline}\label{eq:x_ij_n}x_{ij} = \\\scriptstyle
{\frac{{ \sum\limits_{t = 1}^{p} \sum\limits_{f = 1}^{n}
\sum\limits_{k}n^{(-\frac{1}{2})}_{ik}{{\sum\limits_{\beta \in
J_{r_{1},\,n} {\left\{ {k} \right\}}} {{\rm{cdet}} _{k} \left( {\left(
\widetilde{{\bf A}}^{*}\widetilde{{\bf A}}
 \right)_{\,. \,k}  \left({{\bf e}}_{.f}\right)  } \right)
_{\beta} ^{\beta} } } }
}\,\,\widetilde{d}_{ft}{\sum\limits_{l}{{\sum\limits_{\alpha \in
I_{r_{2},\,p} {\left\{ {l} \right\}}} {{\rm{rdet}} _{l} \left( {\left(
\widetilde{{\bf B}}\widetilde{{\bf B}}^{*}
 \right)_{l. \,}  ({{\bf e}}_{t.})  } \right)
_{\alpha} ^{\alpha} } } }\cdot m^{(\frac{1}{2})}_{lj}} }{{{\sum\limits_{\beta \in
J_{r_{1},\,n}} {{\left| {\left(
\widetilde{{\bf A}}^{*}\widetilde{{\bf A}}
 \right)  _{\beta} ^{\beta} }  \right|}}}
}{{\sum\limits_{\alpha \in I_{r_{2},p}}
{{\left| {\left( \widetilde{{\bf B}}\widetilde{{\bf B}}^{*} \right) _{\alpha} ^{\alpha} }  \right|}}} }}.
 }
\end{multline}
Denote by
\begin{multline*}
 d^{{\rm {\bf A}}}_{it}:=
\sum\limits_{k}n^{(-\frac{1}{2})}_{ik}\sum\limits_{\beta \in
J_{r_{1},\,n} {\left\{ {k} \right\}}} {\rm{cdet}} _{k} \left( {\left(
\widetilde{{\bf A}}^{*}\widetilde{{\bf A}}
 \right)_{\,. \,k} \left( \widetilde{{\rm {\bf d}}}_{.t}
\right)} \right) _{\beta} ^{\beta}  =\\ \sum\limits_{f
= 1}^{n} {\sum\limits_{k}n^{(-\frac{1}{2})}_{ik}\sum\limits_{\beta \in
J_{r_{1},\,n} {\left\{ {k} \right\}}} {\rm{cdet}} _{k} \left( {\left(
\widetilde{{\bf A}}^{*}\widetilde{{\bf A}}
 \right)_{\,. \,k} \left({{\rm {\bf e}}}_{.f}
\right)} \right) _{\beta} ^{\beta}}
\widetilde{d}_{ft}
\end{multline*}
the $t$th component  of the row-vector ${\rm {\bf d}}^{{\rm {\bf
A}}}_{i\,.}= (d^{{\rm {\bf A}}}_{i1},...,d^{{\rm {\bf A}}}_{ip})$
for all $t=1,...,p$. Substituting it in (\ref{eq:x_ij_n}),
we have
\[x_{ij} ={\frac{{ \sum\limits_{t = 1}^{p}
 d^{{\rm {\bf A}}}_{it}
}{\sum\limits_{l}{{\sum\limits_{\alpha \in
I_{r_{2},\,p} {\left\{ {l} \right\}}} {{\rm{rdet}} _{l} \left( {\left(
\widetilde{{\bf B}}\widetilde{{\bf B}}^{*}
 \right)_{l. \,}  ({{\bf e}}_{t.})  } \right)
_{\alpha} ^{\alpha} } } }\cdot m^{(\frac{1}{2})}_{lj}}
}{{{\sum\limits_{\beta \in
J_{r_{1},\,n}} {{\left| {\left(
\widetilde{{\bf A}}^{*}\widetilde{{\bf A}}
 \right)  _{\beta} ^{\beta} }  \right|}}}
}{{\sum\limits_{\alpha \in I_{r_{2},p}} {{\left| {\left( \widetilde{{\bf B}}\widetilde{{\bf B}}^{*} \right) _{\alpha} ^{\alpha} }  \right|}}} }}     }.
\]
Since $\sum\limits_{t = 1}^{p}
 d^{{\rm {\bf A}}}_{it}{\rm {\bf e}}_{t.}={\rm {\bf
d}}^{{\rm {\bf A}}}_{i\,.}$, then it follows (\ref{eq:d^A_n}).

If we denote by
\begin{multline*}
 \sum\limits_{t = 1}^{p}\widetilde{d}_{ft}{\sum\limits_{l}{{\sum\limits_{\alpha \in
I_{r_{2},\,p} {\left\{ {l} \right\}}} {{\rm{rdet}} _{l} \left( {\left(
\widetilde{{\bf B}}\widetilde{{\bf B}}^{*}
 \right)_{l. \,}  ({{\bf e}}_{t.})  } \right)
_{\alpha} ^{\alpha} } } }\cdot m^{(\frac{1}{2})}_{lj}}=\\{\sum\limits_{l}{{\sum\limits_{\alpha \in
I_{r_{2},\,p} {\left\{ {l} \right\}}} {{\rm{rdet}} _{l} \left( {\left(
\widetilde{{\bf B}}\widetilde{{\bf B}}^{*}
 \right)_{l. \,}  ({\widetilde{\bf d}}_{f.})  } \right)
_{\alpha} ^{\alpha} } } }\cdot m^{(\frac{1}{2})}_{lj}}=:d^{{\rm {\bf B}}}_{fj}
\end{multline*}
\noindent the $f$th component  of the column-vector ${\rm {\bf
d}}^{{\rm {\bf B}}}_{.\,j}= (d^{{\rm {\bf B}}}_{1j},...,d^{{\rm
{\bf B}}}_{jn})^{T}$ for all $f=1,...,n$ and substitute it
in (\ref{eq:x_ij_n}), then
\[x_{ij} ={\frac{{   \sum\limits_{f
= 1}^{n} {\sum\limits_{k}n^{(-\frac{1}{2})}_{ik}\sum\limits_{\beta \in
J_{r_{1},\,n} {\left\{ {k} \right\}}} {\rm{cdet}} _{k} \left( {\left(
\widetilde{{\bf A}}^{*}\widetilde{{\bf A}}
 \right)_{\,. \,k} \left({{\rm {\bf e}}}_{.f}
\right)} \right) _{\beta} ^{\beta}} }\,\,d^{{\rm {\bf
B}}}_{fj} }{{{\sum\limits_{\beta \in J_{r_{1},\,n}} {{\left|
{\left( {{\rm {\bf A}}^{\sharp} {\rm {\bf A}}} \right){\kern 1pt}
_{\beta} ^{\beta} }  \right|}}} }{{\sum\limits_{\alpha \in
I_{r_{2},p}} {{\left| {\left( {{\rm {\bf B}}{\rm {\bf B}}^{\sharp} }
\right){\kern 1pt} _{\alpha} ^{\alpha} } \right|}}} }}    }.
\]
Since $\sum\limits_{f = 1}^{n}{\rm {\bf e}}_{.f}
 d^{{\rm {\bf B}}}_{fj}={\rm {\bf
d}}^{{\rm {\bf B}}}_{.\,j}$, then it follows (\ref{eq:d^B_n}).

(ii) If ${\rm rank}\,{\rm {\bf A}} = n$ and ${\rm rank}\,{\rm {\bf
B}} = p$,
  then by Theorem \ref{th:det_rep_A_mn_non} the weighted  Moore-Penrose inverses ${\rm {\bf
A}}_{M,N}^{ \dag} = \left( {{a}_{ij}^{ \ddag} } \right) \in
{\mathbb{H}}^{n\times m} $ and ${\rm {\bf B}}_{P,Q}^{\dag} = \left(
{{b}_{ij}^{\ddag} } \right) \in {\rm {\mathbb{H}}}^{q\times p} $ possess
the following determinantal representations, respectively,
\begin{gather}\label{eq:A1A_f_n}
 {{a}_{ij}^{\ddag} }   = \frac{{\rm{cdet}} _{i}({\bf A}^{*}{\bf M}{\bf A})_{.i}\,(\widehat{\bf{a}}_{.j})}{\det({\bf A}^{*}{\bf M}{\bf A})},
\\\label{eq:BB1_f_n}
  {{b}_{ij}^{\ddag} }   =    \frac{{\rm{rdet}}_{j}({\bf B}{\bf Q}^{-1}{\bf B}^{*})_{j.}(\widehat{\bf{b}}_{i.})}{{\det}({\bf B}{\bf Q}^{-1}{\bf B}^{*})}.
\end{gather}
where $\widehat{\bf{a}}_{.j}$ is the $j$th column of ${\bf A}^{*}{\bf M}$ for all $j=1,\ldots,m$, and $\widehat{\bf{b}}_{i.}$ is the $i$th row of ${\bf Q}^{-1}{\bf B}^{*}$ for all $i=1,\ldots,n$.

By their substituting in (\ref{eq:sum+_n}), we obtain
\begin{equation*}
x_{ij} =
\frac{ \sum\limits_{t = 1}^{p} \sum\limits_{f = 1}^{n}
{\rm{cdet}} _{i}({\bf A}^{*}{\bf M}{\bf A})_{.i}\,({\bf{e}}_{.f})\,\widetilde{d}_{ft}{\rm{rdet}}_{j}({\bf B}{\bf Q}^{-1}{\bf B}^{*})_{j.}({\bf{e}}_{t.}) }{{\rm det}\,({\bf A}^{*}{\bf M}{\bf A}){\det}({\bf B}{\bf Q}^{-1}{\bf B}^{*})},
\end{equation*}
where $\widetilde{d}_{ft}$ is the $(ft)$th entry of $\widetilde{{\bf D}}:={\bf A}^{*}{\bf M}{\bf D}{\bf Q}^{-1}{\bf B}^{*}$ in this case.
Denote by
\begin{equation*}
 d^{{\rm {\bf A}}}_{it}:=
{\rm{cdet}} _{i}({\bf A}^{*}{\bf M}{\bf A})_{.i}\,(\widetilde{\bf{d}}_{.t})
\end{equation*}
the $t$th component  of the row-vector ${\rm {\bf d}}^{{\rm {\bf
A}}}_{i\,.}= (d^{{\rm {\bf A}}}_{i1},...,d^{{\rm {\bf A}}}_{ip})$
for all $t=1,...,p$. Substituting it in (\ref{eq:x_ij_n}),
it follows (\ref{eq:AXB_cdetA*A_n}).

Similarly, we can obtain (\ref{eq:AXB_rdetBB*_n}).

(iii)
 If ${\rm {\bf A}} \in {\rm {\mathbb{H}}}_{r_{1}}^{m\times n}
$, ${\rm {\bf B}} \in {\rm {\mathbb{H}}}_{r_{2}}^{p\times q} $ and
$ r_{1} =n$, $ r_{2} < p$, then, for the weighted  Moore-Penrose inverses ${\bf A}_{M,N}^{\dag}$ and ${\bf B}_{P,Q}^{\dag}$, the determinantal representations (\ref{eq:A1A_f_n})
 and  (\ref{eq:a+nf_n})  are more applicable to use, respectively.
By their substituting in (\ref{eq:sum+_n}) and
pondering ahead as in the previous case, we finally obtain (\ref{eq:AXB_detA*A_d^B_n}) and (\ref{AXB_detA*A_d^A_n}) as well.

(iv) In this case  for  ${\bf A}_{M,N}^{\dag}$ and ${\bf B}_{P,Q}^{\dag}$, we use the determinantal representations
(\ref{eq:a+nf_n}) and (\ref{eq:BB1_f_n}), respectively.
$\Box$

\begin{corollary}\label{cor:AX_n} Suppose that ${\bf A}\in {\mathbb{H}}^{m\times n}_{{r}_1}$, ${\bf D}\in {\mathbb{H}}^{m\times p}$,   ${\bf M}$, ${\bf N}$ are Hermitian positive definite matrices of order $m$  and $n$, respectively, and ${\bf A}^{\sharp}{\bf A}$ is non-Hermitian.
 If ${\bf D}\subset  \mathcal{R}_{r}({\bf A}{\bf A}^{\sharp})$ and ${\bf D}\subset  \mathcal{R}_{l}({\bf A}^{\sharp}{\bf A})$, then the unique solution ${\bf X}= {\bf A}^{\dag}_{M,N}{\bf D}$ of the equation ${\bf A}{\bf X}={\bf D}$
 with the restrictions (\ref{res_AX})  possess the following determinantal representations.
\begin{enumerate}
\item[(i)] If ${\rm rank}\,{\rm {\bf A}} = r_{1} <  n$, then
\begin{equation*}
x_{ij} = {\frac{\sum\limits_{k}n^{(-\frac{1}{2})}_{ik}{{\sum\limits_{\beta \in
J_{r_{1},\,n} {\left\{ {i} \right\}}} {{\rm{cdet}} _{k} \left( {\left(
\widetilde{{\bf A}}^{*}\widetilde{{\bf A}}
 \right)_{\,. \,k}  \left({\widetilde{\bf d}}_{.j}\right)  } \right)
_{\beta} ^{\beta} } } }     }{{{\sum\limits_{\beta \in J_{r_{1},\,n}}
{{\left| {\left(
\widetilde{{\bf A}}^{*}\widetilde{{\bf A}}
 \right)  _{\beta} ^{\beta} }  \right|}}} }}},
\end{equation*}
where
${\widetilde{\bf d}}_{.j}$ are  the $j$th
column  of ${ \widetilde{\bf{D}}}= {\bf N}^{-\frac{1}{2}}{\bf A}^{*}{\bf M}{\bf D}$ for all $i =1,...,n $, $j =1,...,p $.
\item[(ii)] If ${\rm rank}\,{\rm {\bf A}} = n$, then
\begin{equation*}
 x_{i\,j} = \frac{{{\rm cdet} _{i} \left(
{\bf A}^{*}{\bf M}{\bf A}
 \right)_{.\,i\,} \left( {{\widetilde {\bf
d}}_{.j}} \right)}}{{\rm det} \left(
{\bf A}^{*}{\bf M}{\bf A}
 \right)},\end{equation*}
 where
${\widetilde{\bf d}}_{.j}$ are  the $j$th
column  of ${ \widetilde{\bf{D}}}=  {\bf
A}^{*}{\bf M} {\bf D}$.

\end{enumerate}
\end{corollary}
{\bf Proof.} The proof follows evidently from Theorem \ref{th:cr_rul_AXB_n} when ${\bf B}$  be removed and unit matrices  insert instead ${\bf P}$, ${\bf Q}$.
\begin{corollary}\label{cor:XB_n} Suppose that ${\bf B}\in {\mathbb{H}}^{p\times q}_{{r}_2}$, ${\bf D}\in {\mathbb{H}}^{n\times q}$, ${\bf P}$, and ${\bf Q}$ are Hermitian positive definite matrices of order  $p$ and $q$, respectively, and ${\bf B}{\bf B}^{\sharp}$ is non-Hermitian. If ${\bf D}\subset  \mathcal{R}_{r}({\bf B}^{\sharp}{\bf B})$ and ${\bf D}\subset  \mathcal{R}_{l}({\bf B}{\bf B}^{\sharp})$, then the unique solution ${\bf X}= {\bf D}{\bf B}^{\dag}_{P,Q}$ of the equation ${\bf X}{\bf B}={\bf D}$
 with the restrictions (\ref{res_XB})  possess the following determinantal representations.
\begin{enumerate}
\item[(i)] If ${\rm rank}\,{\bf B} = r_{2} <  p$, then
\begin{equation*}
 x_{ij}={\frac{{{\sum\limits_{\alpha
\in I_{r_{2},q} {\left\{ {j} \right\}}} {{\rm{rdet}} _{j} \left(
{\left( {{\widetilde {\bf B}}{\widetilde {\bf B}}^{*} } \right)_{j\,.} \left(
{{{\widetilde {\bf d}}}\,_{i\,.}} \right)}
\right)\,_{\alpha} ^{\alpha} } }}}{{{\sum\limits_{\alpha \in
I_{r_{2},q}} {{\left| {\left( {{\widetilde {\bf B}}{\widetilde {\bf B}}^{*} }
\right) _{\alpha} ^{\alpha} } \right|}}} }}},
\end{equation*}
where
${\widetilde {\bf d}}_{i.}$ are  the $i$th
row  of ${ \widetilde{\bf{D}}}=  {\bf D}{\bf Q}^{-1} {\bf
B}^{*}{\bf P}^{\frac{1}{2}}$  for all $i =1,...,n $, $j =1,...,p $.
\item[(ii)] If ${\rm rank}\,{\bf B} =   p$, then
\begin{equation}\label{eq:XB_nf}
 x_{i\,j} = \frac{{{\rm rdet} _{j} \left({\bf B}{\bf Q}^{-1}{\bf B}^{*}\right)_{j.} \left( {{\widetilde {\bf
d}}_{i.}} \right)}}{{\rm det} \left({\bf B}{\bf Q}^{-1}{\bf B}^{*}\right)},
\end{equation}
 where
${\widetilde{\bf d}}_{i.}$ are  the $i$th
row  of ${ \widetilde{\bf{D}}}= {\bf D}{\bf Q}^{-1} {\bf
B}^{*} $.

\end{enumerate}
\end{corollary}{\bf Proof.} The proof follows evidently from Theorem \ref{th:cr_rul_AXB_n} when ${\bf A}$  be removed and unit matrices  insert instead ${\bf M}$, ${\bf N}$.

\subsection{Mixed Cases}
In this subsection we consider mixed cases when only one from the pair ${\bf A}^{\sharp}{\bf A}$ and ${\bf B}{\bf B}^{\sharp}$ is non-Hermitian.  We give this theorems without proofs, since  their proofs  are similar to the proof of Theorems \ref{th:cr_rul_AXB} and \ref{th:cr_rul_AXB_n}.

\begin{theorem}\label{th:cr_rul_AXB_1n} Let ${\bf A}^{\sharp}{\bf A}$ be Hermitian and ${\bf B}{\bf B}^{\sharp}$ be  non-Hermitian.  Then the  solution  (\ref{eq:sol_res})  possess the following determinantal representations.
\begin{enumerate}
\item[(i)] If ${\rm rank}\,{\rm {\bf A}} = r_{1} <  n$ and
 ${\rm rank}\,{\rm {\bf B}} = r_{2} < p$, then
\begin{equation*}
x_{ij} = {\frac{{{\sum\limits_{\beta \in J_{{{r}_1},\,n} {\left\{
{i} \right\}}} {{\rm{cdet}} _{i} \left( {\left( {{\rm {\bf A}}^{\sharp} {\rm {\bf A}}} \right)_{.\,i} \left( {{{\rm {\bf
d}}}\,_{.\,j}^{{\rm {\bf B}}}} \right)} \right) _{\beta} ^{\beta}
} } }}{{{\sum\limits_{\beta \in J_{{{r}_1},n}} {{\left| {\left(
{{\rm {\bf A}}^{\sharp} {\rm {\bf A}}} \right)_{\beta} ^{\beta} }
\right|}} \sum\limits_{\alpha \in I_{{{r}_2},p}}{{\left| {\left(
\widetilde{{\bf B}}\widetilde{{\bf B}}^{*}
 \right) _{\alpha} ^{\alpha} }
\right|}}} }}},
\end{equation*} or

\begin{equation*}
 x_{ij}={\frac{{\sum\limits_{l}{{\sum\limits_{\alpha \in
I_{r_{2},\,p} {\left\{ {l} \right\}}} {{\rm{rdet}} _{l} \left( {\left(
\widetilde{{\bf B}}\widetilde{{\bf B}}^{*}
 \right)_{l. \,}  ({{\bf d}}^{{\rm {\bf A}}}_{i.})  } \right)
_{\alpha} ^{\alpha} } } }\cdot m^{(\frac{1}{2})}_{lj}}}{{{\sum\limits_{\beta \in
J_{r_{1},n}} {{\left| {\left( {{\rm {\bf A}}^{\sharp} {\rm {\bf A}}}
\right) _{\beta} ^{\beta} } \right|}}\sum\limits_{\alpha \in
I_{r_{2},p}} {{\left| {\left(
\widetilde{{\bf B}}\widetilde{{\bf B}}^{*}
 \right) _{\alpha} ^{\alpha} } \right|}}} }}},
\end{equation*}

where
\begin{equation} \label{eq:def_d^B_m_1n}
   {{{\rm {\bf d}}}_{.\,j}^{{\rm {\bf B}}}}=\left(
{\sum\limits_{l}{{\sum\limits_{\alpha \in
I_{r_{2},\,p} {\left\{ {l} \right\}}} {{\rm{rdet}} _{l} \left( {\left(
\widetilde{{\bf B}}\widetilde{{\bf B}}^{*}
 \right)_{l. \,}  ({\widetilde{\bf d}}_{t.})  } \right)
_{\alpha} ^{\alpha} } } }\cdot m^{(\frac{1}{2})}_{lj}}
\right)\in{\rm {\mathbb{H}}}^{n \times
1}\end{equation}
\begin{equation} \label{eq:d_A<n_1n}  {{{\rm {\bf d}}}_{i\,.}^{{\rm {\bf A}}}}=\left(
\sum\limits_{\beta \in J_{r_{1},n} {\left\{ {i} \right\}}}
{{\rm{cdet}} _{i} \left( {\left( {{\rm {\bf A}}^{\sharp}{\rm {\bf A}}
} \right)_{.i} \left( {\widetilde{{\rm {\bf d}}}_{.f}} \right)}
\right)_{\beta} ^{\beta}} \right)\in{\rm {\mathbb{H}}}^{1 \times
p}
\end{equation}
are the column-vector and the row-vector, respectively.  ${\tilde{{\rm {\bf
 d}}}_{t.}}$ and
${\widetilde{{\rm {\bf d}}}_{.f}}$ are the $t$th row   and the $f$th
column  of ${\widetilde {\bf {D}}}:={{\bf A}}^{\sharp}{\bf D}{\bf Q}^{-1} {\bf
B}^{*}{\bf P}^{\frac{1}{2}}$ for all $t =1,...,n $, $f =1,...,p $.

\item[(ii)] If ${\rm rank}\,{\rm {\bf A}} = n$ and ${\rm rank}\,{\rm {\bf B}} = p$, then
\begin{equation*}
 x_{i\,j} = \frac{{{\rm cdet} _{i} (
{\bf A}^{\sharp}  {\bf A})_{.\,i\,} \left( {{\rm {\bf
d}}_{.j}^{{\rm {\bf B}}}} \right)}}{{\rm det} ({\bf A}^{\sharp}  {\bf A})\cdot {\rm det} ({\bf B}{\bf Q}^{-1}{\bf B}^{*})},
\end{equation*}
or
\begin{equation*}
 x_{i\,j} = \frac{{\rm rdet} _{j} ({\bf B}{\bf Q}^{-1}{\bf B}^{*})_{j.\,} \left(  {\bf
d}_{i\,.}^{ {\bf A}} \right)}{{\rm det} ({\bf A}^{\sharp}  {\bf A})\cdot {\rm det} ({\bf B}{\bf Q}^{-1}{\bf B}^{*})},
\end{equation*}
 \noindent where   \begin{equation}\label{eq:d_B_p_1n} {\bf d}_{.j}^ {\bf B} : =
\left( {\rm
rdet} _{j} ({\bf B}{\bf Q}^{-1}{\bf B}^{*})_{j.\,} \left(
{\tilde {\bf d}}_{t\,.} \right) \right)\in{\mathbb{H}}^{n \times
1},\end{equation}
\begin{equation}\label{eq:def_d^A_1n}{\bf d}_{i\,.}^ {\bf A} : = \left(
{\rm cdet}
_{i} ( {\bf A}^{\sharp}  {\bf A})_{.\,i} \left(
{\tilde {\bf d}}_{.f} \right) \right)\in{\rm {\mathbb{H}}}^{1 \times
p},
\end{equation}
${\tilde {\bf d}}_{t\,.}$, ${\tilde {\bf d}}_{.f}$ are the $t$th row and $f$th column of ${\tilde {\bf D}}={\bf A}^{\sharp}{\bf D}{\bf Q}^{-1}{\bf B}^{*}$.

 \item[(iii)]
If  ${\rm rank}\,{\rm {\bf A}} = n$ and  ${\rm rank}\,{\rm {\bf
B}} = r_{2} < p$, then
\begin{equation*}
 x_{ij}={\frac{{{ {{\rm{cdet}} _{i} \left( {\left( {{\rm {\bf A}}^{\sharp} {\rm {\bf A}}} \right)_{.\,i} \left( {{{\rm {\bf
d}}}\,_{.\,j}^{{\rm {\bf B}}}} \right)} \right)
} } }}{{{\rm det}( {\bf A}^{\sharp}  {\bf A})\cdot\sum\limits_{\alpha \in I_{r_{2},p}}{{\left| {\left( \widetilde{{\bf B}}\widetilde{{\bf B}}^{*} \right) _{\alpha} ^{\alpha} }
\right|}}} }},
\end{equation*}
or
\begin{equation*}
x_{ij}={\frac{{\sum\limits_{l}{{\sum\limits_{\alpha \in
I_{r_{2},\,p} {\left\{ {l} \right\}}} {{\rm{rdet}} _{l} \left( {\left(
\widetilde{{\bf B}}\widetilde{{\bf B}}^{*}
 \right)_{l. \,}  ({{\bf d}}^{{\rm {\bf A}}}_{i.})  } \right)
_{\alpha} ^{\alpha} } } }\cdot m^{(\frac{1}{2})}_{lj}}}{{{\rm det}( {\bf A}^{\sharp}  {\bf A})\cdot\sum\limits_{\alpha \in I_{r_{2},p}}{{\left| {\left( \widetilde{{\bf B}}\widetilde{{\bf B}}^{*} \right) _{\alpha} ^{\alpha} }
\right|}}} }},
\end{equation*}
where  $ {{\rm {\bf d}}_{.\,j}^{{\rm {\bf B}}}}$ is
(\ref{eq:def_d^B_m_1n})
   and
  ${\rm {\bf
d}}_{i\,.}^{{\rm {\bf A}}}$ is (\ref{eq:def_d^A_1n}).

\item[(iv)]
If ${\rm rank}\,{\rm {\bf A}} = r_{1} <  n$ and ${\rm rank}\,{\rm
{\bf B}} =  p$, then
\begin{equation*}
 x_{i\,j} = {\frac{{{\rm rdet} _{j} ({\bf B}{\bf Q}^{-1}{\bf B}^{*})_{j.\,} \left( {{\rm {\bf
d}}_{i\,.}^{{\rm {\bf A}}}} \right)}}{{{\sum\limits_{\beta
\in J_{r_{1},n}} {{\left| {\left( {{\rm {\bf A}}^{\sharp} {\rm {\bf
A}}} \right) _{\beta} ^{\beta} } \right|}}\cdot  {\rm det} ({\bf B}{\bf Q}^{-1}{\bf B}^{*})}}}},
\end{equation*}
or
\begin{equation*} x_{i\,j}=
{\frac{{{\sum\limits_{\beta \in J_{r_{1},\,n} {\left\{
{i} \right\}}} {{\rm{cdet}} _{i} \left( {\left( {{\rm {\bf A}}^{\sharp} {\rm {\bf A}}} \right)_{\,.\,i} \left( {{{\rm {\bf
d}}}\,_{.\,j}^{{\rm {\bf B}}}} \right)} \right) _{\beta} ^{\beta}
} } }}{{{\sum\limits_{\beta
\in J_{r_{1},n}} {{\left| {\left( {{\rm {\bf A}}^{\sharp} {\rm {\bf
A}}} \right)_{\beta} ^{\beta} } \right|}}\cdot {\rm det} ({\bf B}{\bf Q}^{-1}{\bf B}^{*})}}}},
\end{equation*}
 \noindent where  $ {{\rm {\bf d}}_{.\,j}^{{\rm {\bf B}}}}$ is
(\ref{eq:d_B_p_1n})
   and
 ${\rm {\bf
d}}_{i\,.}^{{\rm {\bf A}}}$ is (\ref{eq:d_A<n_1n}).

\end{enumerate}
\end{theorem}

\begin{theorem}\label{th:cr_rul_AXB_n1} Let ${\bf A}^{\sharp}{\bf A}$ be  non-Hermitian, and ${\bf B}{\bf B}^{\sharp}$ be Hermitian. Denote ${\widetilde {\bf {D}}}:=\widetilde{{\bf A}}^{*}{\bf D}{\bf B}^{\sharp}$. Then the  solution  (\ref{eq:sol_res})  possess the following determinantal representations.
\begin{enumerate}
\item[(i)] If ${\rm rank}\,{\rm {\bf A}} = r_{1} <  n$ and
 ${\rm rank}\,{\rm {\bf B}} = r_{2} < p$, then
\begin{equation*}
x_{ij} = {\frac{{    {\sum\limits_{k}n^{(-\frac{1}{2})}_{ik}\sum\limits_{\beta \in
J_{r_{1},\,n} {\left\{ {k} \right\}}} {\rm{cdet}} _{k} \left( {\left(
\widetilde{{\bf A}}^{*}\widetilde{{\bf A}}
 \right)_{\,. \,k} \left({{\rm {\bf d}}}^{{\rm {\bf
B}}}_{.j}
\right)} \right) _{\beta} ^{\beta}} } }{{{\sum\limits_{\beta \in J_{r_{1},\,n}} {{\left|
{\left(
\widetilde{{\bf A}}^{*}\widetilde{{\bf A}}
 \right){\kern 1pt}
_{\beta} ^{\beta} }  \right|}}} }{{\sum\limits_{\alpha \in
I_{r_{2},p}} {{\left| {\left(
{{\rm {\bf B}}{\rm {\bf B}}^{\sharp} } \right) _{\alpha} ^{\alpha} } \right|}}} }}    },
\end{equation*} or

\begin{equation*}
 x_{ij}={\frac{{{\sum\limits_{\alpha
\in I_{r_{2},p} {\left\{ {j} \right\}}} {{\rm{rdet}} _{j} \left(
{\left( {{\rm {\bf B}}{\rm {\bf B}}^{\sharp} } \right)_{\,j\,.} \left(
{{{\rm {\bf d}}}\,_{i\,.}^{{\rm {\bf A}}}} \right)}
\right)\,_{\alpha} ^{\alpha} } }}
}{{{\sum\limits_{\beta \in
J_{r_{1},\,n}} {{\left| {\left(
\widetilde{{\bf A}}^{*}\widetilde{{\bf A}}
 \right)  _{\beta} ^{\beta} }  \right|}}}
}{{\sum\limits_{\alpha \in I_{r_{2},p}} {{\left| {\left( {{\rm {\bf B}}{\rm {\bf B}}^{\sharp} } \right) _{\alpha} ^{\alpha} }  \right|}}} }}     },
\end{equation*}

where
\begin{equation} \label{eq:def_d^B_m_n1}
   {{{\rm {\bf d}}}_{.\,j}^{{\rm {\bf B}}}}=\left(
{{\sum\limits_{\alpha
\in I_{r_{2},p} {\left\{ {j} \right\}}} {{\rm{rdet}} _{j} \left(
{\left( {{\rm {\bf B}}{\rm {\bf B}}^{\sharp} } \right)_{\,j\,.} \left(
{{{\widetilde {\bf d}}}\,_{t\,.}
} \right)}
\right)\,_{\alpha} ^{\alpha} } }}
\right)\in{\rm {\mathbb{H}}}^{n \times
1}\end{equation}
\begin{equation} \label{eq:d_A<n_n1}  {{{\rm {\bf d}}}_{i\,.}^{{\rm {\bf A}}}}=\left(
 {\sum\limits_{k}n^{(-\frac{1}{2})}_{ik}\sum\limits_{\beta \in
J_{r_{1},\,n} {\left\{ {k} \right\}}} {\rm{cdet}} _{k} \left( {\left(
\widetilde{{\bf A}}^{*}\widetilde{{\bf A}}
 \right)_{\,. \,k} \left({{\widetilde {\bf d}}}_{.f}
\right)} \right) _{\beta} ^{\beta}}
 \right)\in{\rm {\mathbb{H}}}^{1 \times
p}
\end{equation}
are the column-vector and the row-vector, respectively.  ${\tilde{{\rm {\bf
 d}}}_{t.}}$ and
${\tilde{{\rm {\bf d}}}_{.f}}$ are the $t$th row   and the $f$th
column  of ${\widetilde {\bf {D}}}:={\bf N}^{-\frac{1}{2}}{\bf A}^{*}{\bf M}{\bf D}{\bf B}^{\sharp}$ for all $t =1,...,n $, $f =1,...,p $.

\item[(ii)] If ${\rm rank}\,{\rm {\bf A}} = n$ and ${\rm rank}\,{\rm {\bf B}} = p$, then
\begin{equation*}
 x_{i\,j} = \frac{{{\rm cdet} _{i} ({\bf A}^{*}{\bf M}{\bf A})_{.\,i\,} \left( {{\rm {\bf
d}}_{.j}^{{\rm {\bf B}}}} \right)}}{{\rm det} ({\bf A}^{*}{\bf M}{\bf A})\cdot {\rm det} \left( {{\rm {\bf B}}{\rm {\bf B}}^{\sharp} } \right)},
\end{equation*}
or
\begin{equation*}
 x_{i\,j} = \frac{{\rm rdet} _{j} \left( {{\rm {\bf B}}{\rm {\bf B}}^{\sharp} } \right)_{j.\,} \left(  {\bf
d}_{i\,.}^{ {\bf A}} \right)}{{\rm det} ({\bf A}^{*}{\bf M}{\bf A})\cdot {\rm det} \left( {{\rm {\bf B}}{\rm {\bf B}}^{\sharp} } \right)},
\end{equation*}
 \noindent where   \begin{equation}\label{eq:d_B_p_n1} {\bf d}_{.j}^ {\bf B} : =
\left( {\rm
rdet} _{j} \left( {{\rm {\bf B}}{\rm {\bf B}}^{\sharp} } \right)_{j.\,} \left(
{\tilde {\bf d}}_{t\,.} \right) \right)\in{\mathbb{H}}^{n \times
1},\end{equation}
\begin{equation}\label{eq:def_d^A_n1}{\bf d}_{i\,.}^ {\bf A} : = \left(
{\rm cdet}
_{i} ({\bf A}^{*}{\bf M}{\bf A})_{.\,i} \left(
{\tilde {\bf d}}_{.f} \right) \right)\in{\rm {\mathbb{H}}}^{1 \times
p},\end{equation}
${\tilde {\bf d}}_{t\,.}$, ${\tilde {\bf d}}_{.f}$ are the $t$th row and $f$th column of ${\tilde {\bf D}}={\bf A}^{*}{\bf M}{\bf D}{\bf B}^{\sharp}$.
 \item[(iii)]
If  ${\rm rank}\,{\rm {\bf A}} = n$ and  ${\rm rank}\,{\rm {\bf
B}} = r_{2} < p$, then
\begin{equation*}
 x_{ij}={\frac{{{ {{\rm{cdet}} _{i} \left( {({\bf A}^{*}{\bf M}{\bf A})_{.\,i} \left( {{{\rm {\bf
d}}}\,_{.\,j}^{{\rm {\bf B}}}} \right)} \right)
} } }}{{{\rm det}({\bf A}^{*}{\bf M}{\bf A})\cdot\sum\limits_{\alpha \in I_{r_{2},p}}{{\left| {\left( {{\rm {\bf B}}{\rm {\bf B}}^{\sharp} } \right) _{\alpha} ^{\alpha} }
\right|}}} }},
\end{equation*}
or
\begin{equation*}
x_{ij}={\frac{{{\sum\limits_{\alpha \in I_{r_{2},p} {\left\{ {j}
\right\}}} {\rm rdet}_{j}{ \left( {\left( {{\rm {\bf B}}{\rm {\bf B}}^{\sharp} } \right)_{\,j\,.} \left( {{{\rm {\bf d}}}\,_{i\,.}^{{\rm
{\bf A}}}} \right)} \right)\,_{\alpha} ^{\alpha} } }}}{{{\rm det}({\bf A}^{*}{\bf M}{\bf A})\cdot\sum\limits_{\alpha \in I_{r_{2},p}}{{\left| {\left(
{{\rm {\bf B}}{\rm {\bf B}}^{\sharp} } \right) _{\alpha} ^{\alpha} }
\right|}}} }},
\end{equation*}
where  $ {{\rm {\bf d}}_{.\,j}^{{\rm {\bf B}}}}$ is
(\ref{eq:def_d^B_m_n1})
   and
  ${\rm {\bf
d}}_{i\,.}^{{\rm {\bf A}}}$ is (\ref{eq:def_d^A_n1}).

\item[(iv)]
If ${\rm rank}\,{\rm {\bf A}} = r_{1} <  n$ and ${\rm rank}\,{\rm
{\bf B}} =  p$, then
\begin{equation*}
 x_{i\,j} = {\frac{{{\rm rdet} _{j} \left( {{\rm {\bf B}}{\rm {\bf B}}^{\sharp} } \right)_{j.\,} \left( {{\rm {\bf
d}}_{i\,.}^{{\rm {\bf A}}}} \right)}}{{{\sum\limits_{\beta
\in J_{r_{1},n}} {{\left| {\left(
\widetilde{{\bf A}}^{*}\widetilde{{\bf A}}
 \right) _{\beta} ^{\beta} } \right|}}\cdot  {\rm det} \left( {{\rm {\bf B}}{\rm {\bf B}}^{\sharp} } \right)}}}},
\end{equation*}
or
\begin{equation*}  x_{i\,j}=
{\frac{{{\sum\limits_{\beta \in J_{{{r}_1},\,n} {\left\{
{i} \right\}}} {{\rm{cdet}} _{i} \left( {\left( {{\rm {\bf A}}^{\sharp} {\rm {\bf A}}} \right)_{.\,i} \left( {{{\rm {\bf
d}}}\,_{.\,j}^{{\rm {\bf B}}}} \right)} \right) _{\beta} ^{\beta}
} } }}{{{\sum\limits_{\beta
\in J_{r_{1},n}} {{\left| {\left(
\widetilde{{\bf A}}^{*}\widetilde{{\bf A}}
 \right)_{\beta} ^{\beta} } \right|}}\cdot {\rm det} ({\bf B}{\bf Q}^{-1}{\bf B}^{*})}}}},
\end{equation*}
 \noindent where  $ {{\rm {\bf d}}_{.\,j}^{{\rm {\bf B}}}}$ is
(\ref{eq:d_B_p_n1})
   and
 ${\rm {\bf
d}}_{i\,.}^{{\rm {\bf A}}}$ is (\ref{eq:d_A<n_n1}).

\end{enumerate}
\end{theorem}

\subsection{An example}
Let
us consider the restricted matrix equation
\begin{equation}\label{eq:ex}
  \begin{gathered}
{\bf A}{\bf X}{\bf B}={\bf D},\\
\mathcal{R}_{r}({\bf X})\subset {\bf N}^{-1}\mathcal{R}_{r}({\bf A}^{*}),\,\mathcal{N}_{r}({\bf X})\supset {\bf P}^{-1}\mathcal{N}_{r}({\bf B}^{*}),\\\mathcal{R}_{l}({\bf X})\subset \mathcal{R}_{l}({\bf A}^{*}){\bf M},\,\mathcal{N}_{l}({\bf X})\supset \mathcal{N}_{l}({\bf B}^{*}){\bf Q}
\end{gathered}
\end{equation}
where \[{\bf A}=\begin{pmatrix}
  k & -j & j \\
1 & 0 & k
\end{pmatrix},\,
{\bf N}=\begin{pmatrix}
  5 & 0 & -4j \\
  0 & 4 & 0 \\
  4j & 0 & 5
\end{pmatrix},\,  {\bf M}=\begin{pmatrix}
 5 & 4k  \\
  -4k & 5  \\
\end{pmatrix},\, {\bf D}=\begin{pmatrix}
  i & -j & k \\
-k & 0 & j
\end{pmatrix},\]
\[{\bf B}=\begin{pmatrix}
  k & -j & j \\
0 & 1 & i
\end{pmatrix},\
{\bf Q}=\begin{pmatrix}
  1 & i & 0 \\
  -i & 2 & -j \\
  0 & j & 2
\end{pmatrix},\,\,  {\bf P}=\begin{pmatrix}
  2.5 & -1.5j  \\
  1.5j & 2.5  \\
\end{pmatrix}.\]
Since
\[{\bf A}^{*}=\begin{pmatrix}
  -k & 1 \\
  i &  0 \\
  -j & -k
\end{pmatrix},\,\, {\bf B}^{*}=\begin{pmatrix}
  -k & 0 \\
  j &  1 \\
  -j & -i
\end{pmatrix},\] \[\det({\bf A}{\bf A}^{*})=\det \begin{pmatrix}
  3 & -i+k \\
  i-k &  2 \\
\end{pmatrix}=4,\\\det({\bf B}{\bf B}^{*})=\det \begin{pmatrix}
  3 & -j+k \\
  j-k &  2 \\
\end{pmatrix}=4,\] then ${\rm rank}{\bf A}={\rm rank}{\bf B}=2$.

Due to Theorem \ref{kyrc10}, we can be obtain the inverses
\[
{\bf Q}^{-1}=\frac{1}{3}\begin{pmatrix}
  5 & -4i & -3k \\
  4i & 5 & 3j \\
  3k & -3j & 3
\end{pmatrix},\,\,{\bf N}^{-1}=\frac{1}{9}\begin{pmatrix}
  \frac{5}{9} & 0 & \frac{4j}{9} \\
  0 & \frac{1}{4} & 0 \\
  -\frac{4j}{9} & 0 & \frac{5}{9}
\end{pmatrix}.\]
It is easy to verify that the both matrices ${\bf A}^{\sharp}{\bf A}={\bf N}^{-1}{\bf A}^{*}{\bf M}{\bf A}$ and
${\bf B}{\bf B}^{\sharp}={\bf B}{\bf Q}^{-1}{\bf B}^{*}{\bf P}$
are not Hermitian. Hence,
 we shall find the solution of (\ref{eq:ex}) by (\ref{eq:AXB_detBB*_d^A_n}). So,
  \begin{equation}
\label{eq:AXB_detBB*_d^A_n_ex} x_{i\,j} = {\frac{{{\rm rdet} _{j} ({\bf B}{\bf Q}^{-1}{\bf B}^{*})_{j.\,} \left( {{\rm {\bf
d}}_{i\,.}^{{\rm {\bf A}}}} \right)}}{{{\sum\limits_{\beta
\in J_{2,3}} {{\left| {\left(
\widetilde{{\bf A}}^{*}\widetilde{{\bf A}}
 \right) _{\beta} ^{\beta} } \right|}}\cdot  {\rm det} ({\bf B}{\bf Q}^{-1}{\bf B}^{*})}}}},\, i=1,2,3,\,\,j=1,2,
\end{equation}
 where ${\rm {\bf
d}}_{i\,.}^{{\rm {\bf A}}}$ is (\ref{eq:d_A<n_n}), namely
 \begin{equation} \label{eq:d_A<n_n_ex}  {{{\rm {\bf d}}}_{i\,.}^{{\rm {\bf A}}}}=\left(
 {\sum\limits_{k}n^{(-\frac{1}{2})}_{ik}\sum\limits_{\beta \in
J_{2,\,3} {\left\{ {k} \right\}}} {\rm{cdet}} _{k} \left( {\left(
\widetilde{{\bf A}}^{*}\widetilde{{\bf A}}
 \right)_{\,. \,k} \left({{\widetilde {\bf d}}}_{.f}
\right)} \right) _{\beta} ^{\beta}}
 \right)\in{\rm {\mathbb{H}}}^{1 \times
2}.
\end{equation}
 To obtain ${\bf N}^{-\frac{1}{2}}$, we firstly find the eigenvalues of ${\bf N}$ which are the roots of the characteristic
polynomial
\[
p(\lambda)=\det\begin{pmatrix}\lambda-5 & 0 & 4j\\
 0 & \lambda-4 & 0\\
                                                      -4j & 0&  \lambda-5
\end{pmatrix}=\lambda^{3}-14\lambda^{2}+49\lambda-36\,\,\Rightarrow\,\,\begin{cases}
                                                           \lambda_{1}=1, \\
                                                           \lambda_{2}=4,\\
                                                           \lambda_{3}=9.
                                                         \end{cases}
\]
 By computing the associated eigenvectors and after their orthonormalization, we obtain the unitary matrix ${\bf U}$ whose columns are this eigenvectors.
\[{\bf U}=\begin{pmatrix}0.5 -0.5j & 0 & 0.5 +0.5j\\
 0 & 1 & 0\\
                                                      0.5 +0.5j & 0&  0.5 -0.5j
\end{pmatrix}.\] Finally, we have $${\bf N}^{\frac{1}{2}}={\bf U}^{*}{\bf D}{\bf U}=\begin{pmatrix}
  2 & 0 & -j \\
  0 & 2 & 0 \\
  j & 0 & 2
\end{pmatrix},$$ where ${\bf D}={\rm diag}(1,2,3)$ is the diagonal matrix with $1, 2, 3$ on the principal diagonal. Then by Theorem \ref{kyrc10},
 $${\bf N}^{-\frac{1}{2}}=\begin{pmatrix}
  \frac{2}{3} & 0 & \frac{j}{3} \\
  0 & \frac{1}{2} & 0 \\
   -\frac{j}{3} & 0 & \frac{2}{3}
\end{pmatrix}.$$
 Similarly, ${\bf M}^{\frac{1}{2}}=\begin{pmatrix}
 2 & k  \\
  -k & 2  \\
\end{pmatrix},\,{\bf P}^{\frac{1}{2}}=\frac{1}{2}\begin{pmatrix}
 3 &  -j \\
    j &  3
\end{pmatrix}.$
  Further, we find
 \begin{multline*}\label{Dt}
  \widetilde{{\bf D}}:={\bf N}^{-\frac{1}{2}}{\bf A}^{*}{\bf M}{\bf D}{\bf Q}^{-1}{\bf B}^{*}{\bf P}^{\frac{1}{2}}=\\\begin{pmatrix}
  -\frac{7}{6}+ \frac{125}{3}i+\frac{143}{6}j-\frac{187}{3}k & \frac{47}{2}- \frac{133}{3}i+\frac{61}{6}j-\frac{67}{3}k \\
    \frac{15}{2}+ \frac{43}{2}i+\frac{13}{2}j-\frac{33}{2}k & \frac{15}{2}+ \frac{19}{2}i-\frac{17}{2}j-\frac{29}{2}k \\
   \frac{35}{3}+ \frac{289}{6}i-\frac{101}{3}j+\frac{149}{6}k  & -\frac{47}{3}+ \frac{79}{6}i-11j-\frac{235}{6}k
\end{pmatrix},\\
  \widetilde{{\bf A}}={\bf M}^{\frac{1}{2}}{\bf A}{\bf N}^{-\frac{1}{2}}=\begin{pmatrix}
  \frac{2}{3}+\frac{1}{3}j+2k & -i & -\frac{2}{3}-i+\frac{4}{3}k \\
  2+\frac{2}{3}i-\frac{1}{3}k &\frac{1}{2}j &  \frac{2}{3}i+j+\frac{4}{3}k \\
\end{pmatrix},\\ \widetilde{{\bf A}}^{*}{\bf A}=\begin{pmatrix}
  \frac{82}{9} & -\frac{5}{6}i+3j-\frac{2}{3}k & 3i+\frac{56}{9}j+3k \\
   \frac{5}{6}i-3j+\frac{2}{3}k & \frac{5}{4}& \frac{3}{2}-\frac{4}{3}i+\frac{5}{3}k  \\
   -3i-\frac{56}{9}j-3k & \frac{3}{2}+\frac{4}{3}i-\frac{5}{3}k & \frac{58}{9}
\end{pmatrix}.
\end{multline*}
\begin{multline*}\sum\limits_{\beta
\in J_{2,3}}\left| {\left(
\widetilde{{\bf A}}^{*}\widetilde{{\bf A}}
 \right) _{\beta} ^{\beta} } \right|= \det\begin{pmatrix}\frac{82}{9} & -\frac{5}{6}i+3j-\frac{2}{3}k\\
\frac{5}{6}i-3j+\frac{2}{3}k & \frac{5}{4}\end{pmatrix}+\\\det\begin{pmatrix}\frac{5}{4}& \frac{3}{2}-\frac{4}{3}i+\frac{5}{3}k\\
\frac{3}{2}+\frac{4}{3}i-\frac{5}{3}k & \frac{58}{9}\end{pmatrix}+\det\begin{pmatrix}\frac{82}{9}& 3i+\frac{56}{9}j+3k\\
-3i-\frac{56}{9}j-3k& \frac{58}{9}\end{pmatrix}=\\\frac{5}{4}+\frac{5}{4}+2=4.5,\\
\det \left({\bf B}{\bf Q}^{-1}{\bf B}^{*}\right)= \det \begin{pmatrix}10& 1-2i-4j+k\\
1+2i+4j-k & 3\end{pmatrix}=8.
\end{multline*}
 Now, we obtain components of the row-vectors (\ref{eq:d_A<n_n_ex}),  ${{{\rm {\bf d}}}_{1\,.}^{{\rm {\bf A}}}}=\left(d^{\bf A}_{11}, d^{\bf A}_{12}\right)$,
 \begin{multline*}
 d^{\bf A}_{11}=\frac{2}{3}\left({\rm cdet}_{1}\begin{pmatrix}-\frac{7}{6}+ \frac{125}{3}i+\frac{143}{6}j-\frac{187}{3}k & -\frac{5}{6}i+3j-\frac{2}{3}k\\
\frac{15}{2}+ \frac{43}{2}i+\frac{13}{2}j-\frac{33}{2}k  & \frac{5}{4}\end{pmatrix}+\right.\\\left.{\rm cdet}_{1}\begin{pmatrix}-\frac{7}{6}+ \frac{125}{3}i+\frac{143}{6}j-\frac{187}{3}k & 3i+\frac{56}{9}j+3k\\
\frac{35}{3}+ \frac{289}{6}i-\frac{101}{3}j+\frac{149}{6}k& \frac{58}{9}\end{pmatrix}\right)-\\\frac{1}{3}j\left({\rm cdet}_{2}
\begin{pmatrix}\frac{82}{9}& -\frac{7}{6}+ \frac{125}{3}i+\frac{143}{6}j-\frac{187}{3}k\\
-3i-\frac{56}{9}j-3k& \frac{35}{3}+ \frac{289}{6}i-\frac{101}{3}j+\frac{149}{6}k\end{pmatrix}\right.+\\
\left.{\rm cdet}_{2}\begin{pmatrix}\frac{5}{4}& \frac{15}{2}+ \frac{43}{2}i+\frac{13}{2}j-\frac{33}{2}k \\
\frac{3}{2}+\frac{4}{3}i-\frac{5}{3}k & \frac{35}{3}+ \frac{289}{6}i-\frac{101}{3}j+\frac{149}{6}k\end{pmatrix}\right)=
-\frac{55}{12}- \frac{587}{24}i+\frac{13}{6}j-\frac{253}{8}k.
\end{multline*}
Similarly, we obtain $ d^{\bf A}_{12}=-\frac{13}{12}- \frac{581}{24}i+\frac{39}{4}j+\frac{161}{8}k$. Moreover,
 \begin{multline*}
 {{{\rm {\bf d}}}_{2\,.}^{{\rm {\bf A}}}}=\left(\frac{232}{3}+ \frac{1763}{9}i+\frac{1093}{18}j+\frac{1795}{12}k,
  -\frac{419}{6}+ \frac{3035}{36}i-\frac{751}{9}j-\frac{1189}{9}k\right)\\
  {{{\rm {\bf d}}}_{3\,.}^{{\rm {\bf A}}}}=\left(\frac{116}{9}- \frac{239}{12}i+\frac{25}{36}j+\frac{299}{9}k,
 - \frac{117}{24}+ \frac{95}{6}i-\frac{99}{8}j+\frac{245}{12}k\right).
 \end{multline*}
Finally, we have

\begin{multline*}
 x_{11} = {\frac{{{\rm rdet} _{1} ({\bf B}{\bf Q}^{-1}{\bf B}^{*})_{1.\,} \left( {{\rm {\bf
d}}_{1\,.}^{{\rm {\bf A}}}} \right)}}{{{\sum\limits_{\beta
\in J_{2,3}} {{\left| {\left(
\widetilde{{\bf A}}^{*}\widetilde{{\bf A}}
 \right) _{\beta} ^{\beta} } \right|}}\cdot  {\rm det} ({\bf B}{\bf Q}^{-1}{\bf B}^{*})}}}}=\\\frac{1}{36}\,
 {\rm rdet} _{1}\begin{pmatrix}-\frac{55}{12}- \frac{587}{24}i+\frac{13}{6}j-\frac{253}{8}k& -\frac{13}{12}- \frac{581}{24}i+\frac{39}{4}j+\frac{161}{8}k\\
1+2i+4j-k & 3\end{pmatrix}=\\-\frac{1013}{864}+ \frac{1}{144}i-\frac{359}{864}j+\frac{173}{144}k.
\end{multline*}

Similarly, we obtain
\[x_{12} = \frac{19}{288}- \frac{2459}{864}i-\frac{257}{864}j+\frac{3247}{864}k,\]
and
 \begin{align*}
 {x}_{21}=& \frac{1162}{324}- \frac{8935}{1296}i-\frac{5983}{1296}j+\frac{1759}{432}k,&
{x}_{22}=&  \frac{1631}{432}+ \frac{1285}{324}i-\frac{817}{324}j-\frac{10631}{432}k,\\
 {x}_{31}=&\frac{127}{864}+ \frac{83}{864}i-\frac{545}{864}j-\frac{329}{864}k,&
 {x}_{32}=& \frac{311}{1296}+ \frac{367}{162}i-\frac{949}{1296}j+\frac{77}{36}k.
\end{align*}

Note that we used  Maple with the package CLIFFORD in the calculations.

\section{Conclusion} In this paper, previously obtained determinantal representations  of the quaternion weighted Moore-Penrose inverse have been used to derive explicit determinantal representation formulas  for the solution of the two-sided restricted quaternionic matrix equation, ${\bf A}{\bf X}{\bf B}={\bf D}$, within the framework of the theory of column-row determinants (also previously  introduced by the author).


\begin{thebibliography}{40}\bibitem{pen}  R. A. Penrose, Generalized inverse
 for matrices, {Proc. Camb. Philos. Soc.} \textbf{51} (1955)  406--413. 


\bibitem{ben}
 A. Ben-Israel, T.N.E. Grenville,   {Generalized Inverses: Theory and Applications.} Springer-
Verlag, Berlin, 2002.
\bibitem{pr} K.M. Prasad,  R.B. Bapat, A note of the Khatri inverse, {Sankhya: Indian J. Stat.} \textbf{54} (1992)
291-295. 

    \bibitem{zha}
F. Zhang,   Quaternions and matrices of quaternions, {Linear
Algebra Appl.} \textbf{251} (1997) 21-57. 

   \bibitem{ky_math_sci} I.I. Kyrchei,  Determinantal representation of the Moore-Penrose inverse matrix over the quaternion skew field, {J. Math. Sci.} \textbf{180}(1) (2012) 23-33. 

\bibitem{loan}
C.F. Van Loan,   Generalizing the singular value decomposition, \emph{SIAM J. Numer. Anal.} \textbf{13} (1976)  76--83.

\bibitem{galba}
  E. F. Galba, Weighted singular decomposition and weighted pseudoinversion of matrices, \emph{Ukr. Math. J.} \textbf{48}(10) (1996) 1618-1622. 

\bibitem{kyr_wmpi} I.I. Kyrchei,
 Weighted singular value decomposition and determinantal representations of the quaternion weighted Moore-Penrose inverse, {Appl. Math. Comput.} \textbf{309} (2017)   1-16.

\bibitem{st1}  P. Stanimirovic', M. Stankovic', Determinantal representation of
weighted Moore-Penrose inverse, {Mat. Vesnik}  \textbf{46} (1994) 41-50.

\bibitem{liu2}   X. Liu, Y. Yu,   H. Wang, Determinantal representation of weighted generalized inverses, {Appl. Math. Comput.} \textbf{218}(7) (2011) 3110-3121.

\bibitem{liu1}X. Liu, G. Zhu,   G. Zhou, Y. Yu,  An analog of the adjugate matrix for
the outer inverse ${\bf A}_{T,S}^{(2)}$, {Math. Problem. in Eng. } \textbf{2012},   Article ID 591256 (2012) 14 pages. 

\bibitem{ky_lma1} I.I. Kyrchei,
Analogs of the adjoint matrix for generalized inverses and corresponding
Cramer rules, {Linear Multilinear Algebra} \textbf{56} (4) (2008) 453-469.

 \bibitem{ky_nova} I.I. Kyrchei, Cramer's rule for generalized inverse solutions, In: {Advances in Linear Algebra Research}, I.I.  Kyrchei (Ed.), Nova Sci. Publ., New York,  pp. 79-132, 2015.










 \bibitem{kyr2} I.I.~Kyrchei, {Cramer's rule for quaternion
 systems of linear equations}. Fundamentalnaya i Prikladnaya Matematika \textbf{13}(4) (2007) 67-94.


  \bibitem{ky_th}  I.I. Kyrchei,  The theory of the column and row determinants in a quaternion linear algebra, In: Advances in Mathematics Research \textbf{15}, A.R. Baswell (Ed.), Nova Sci. Publ., New York, pp. 301-359, 2012.



 \bibitem{song1}G. Song,  Q. Wang, H. Chang,   Cramer rule for the unique solution of restricted matrix equations over the quaternion skew field, {Comput Math. Appl.} \textbf{61} (2011)  1576-1589. 

\bibitem{song6}
 G.J. Song,  Determinantal representation of the generalized inverses over
the quaternion skew field with applications, {Appl. Math. Comput.} \textbf{219} (2012)  656--667.


 \bibitem{kyr5}     I.I. Kyrchei,
Explicit representation formulas for the minimum norm least squares solutions of some quaternion matrix equations, {Linear Algebra Appl. } \textbf{438}(1) (2013) 136--152.

 \bibitem{kyr6} I.I. Kyrchei,
Determinantal representations of the Drazin inverse over the quaternion skew field with applications to some matrix equations, {Appl. Math. Comput.} \textbf{238} (2014)  193--207.

 \bibitem{kyr7}I.I. Kyrchei,
Determinantal representations of the W-weighted Drazin inverse over the quaternion skew field, {Appl. Math. Comput.} \textbf{264} (2015)  453--465.

\bibitem{kyr8} I.I. Kyrchei,
Explicit determinantal representation formulas of W-weighted Drazin inverse solutions of some matrix equations over the quaternion skew field, {Math. Problem. in Eng.}  \textbf{2016}  Article ID 8673809 (2016) 13 pages.




\bibitem{kyr9}I.I. Kyrchei,
 Determinantal representations of the Drazin and W-weighted Drazin inverses over the quaternion skew field with applications, In: {Quaternions: Theory and Applications},  S. Griffin (Ed.), New York: Nova Sci. Publ., pp.201-275, 2017.


\bibitem{ky_nova1} A. Kleyn,  I. Kyrchei,  Relation of row-column determinants with quasideterminants of matrices over a quaternion algebra, In: {Advances in Linear Algebra Research},I. Kyrchei (Ed.), Nova Sci. Publ., New York, pp. 299-324, 2015.

\bibitem{song3}G.J. Song,  C.Z.  Dong, New results on condensed Cramer's rule for the general solution to some restricted quaternion matrix equations, {J. Appl. Math. Comput.} \textbf{53}  (2017) 321--341.

  \bibitem{song4}G.J. Song,   Bott-Duffin inverse over the quaternion skew field with applications,
{J. Appl. Math. Comput.} \textbf{41} (2013)  377-392.
 \bibitem{song5}
G.J. Song,  Characterization of the W-weighted Drazin inverse over the quaternion
skew field with applications, {Electron. J. Linear Algebra} \textbf{26} (2013)  1--14.






\bibitem{wei_rep}
Y. Wei,  H.  Wu,  The representation and approximation for the weighted Moore-Penrose inverse, {Appl. Math. Comput.} \textbf{121} (2001)  17-28.
\bibitem{serg}
 I.V. Sergienko,  E.F. Galba,   V.S.  Deineka,Limiting representations of weighted
pseudoinverse matrices with positive
definite weights. Problem regularization, {Cybernetics and Systems Analysis} \textbf{39}(6)  (2003) 816-830.








\bibitem{hu}
 L. Huang, W. So,  On left eigenvalues of a quaternionic
matrix, {Linear Algebra Appl.} \textbf{323}  (2001) 105-116.
 \bibitem{so}
 W. So,  Quaternionic left eigenvalue problem, {Southeast
Asian Bulletin of Mathematics} \textbf{29} (2005)  555-565.
 \bibitem{wo}
R. M. W. Wood,  Quaternionic eigenvalues, {Bull. Lond. Math. Soc.} \textbf{17} (1985)  137-138.





\bibitem{br}
J.L. Brenner,   Matrices of quaternions, {Pac. J. Math. } \textbf{1} (1951) 329-335.
 \bibitem{ma} E. Mac\'{i}as-Virg\'{o}s,  M.J. Pereira-S\'{a}ez,
A topological approach to left eigenvalues of quaternionic matrices,
{Linear Multilinear Algebra}\textbf{ 62}(2)  (2014) 139--158.
\bibitem{ba}A. Baker,  Right eigenvalues for quaternionic matrices: a topological ap-
proach, {Linear Algebra Appl.} \textbf{286} (1999) 303-309.
 \bibitem{dra}T. Dray,  C. A. Manogue,   The octonionic eigenvalue problem, {Advances
in Applied Clifford Algebras} \textbf{8}(2) (1998) 341-364.
\bibitem{zh}F. Zhang,  Quaternions and matrices of quaternions, {Linear Algebra
Appl.} \textbf{251} (1997)  21-57.
\bibitem{far}D.R. Farenick,   B.A.F. Pidkowich, The spectral theorem in quaternions, {Linear
Algebra Appl.} \textbf{ 371} (2003)  75-102.
\bibitem{fa}F. O. Farid,   Q.W.  Wang,   F. Zhang,  On the eigenvalues of quaternion matrices,
{Linear Multilinear Algebra} \textbf{59}(4) (2011)  451--473.



 \bibitem{ho}R.A. Horn, C.R. Johnson,
{Matrix Analysis}. Cambridge etc., Cambridge University
Press 1985.







\end{thebibliography}
\end{document}